\documentclass[12pt]{iopart}

\usepackage{amstext}
\usepackage{bm}
\usepackage{graphicx}
\usepackage{subfig}
\usepackage{algorithm,algorithmic}
\usepackage{tabularx}
\usepackage{tikz, pgfplotstable}

\usetikzlibrary{calc,trees,positioning,arrows,chains,shapes.geometric,decorations.pathreplacing,decorations.pathmorphing,shapes,matrix,shapes.symbols,fillbetween, spy}

\pgfplotsset{every axis/.append style={
		label style={font=\footnotesize},
		ticklabel style={font=\footnotesize, scaled ticks = false},
}}

\definecolor{lightblue}{rgb}{0.68, 0.85, 0.9}

\usepackage{iopams} 
 
\begin{document}
\title[HODLR approximation of Hessians in ice sheet inverse
  problems]{Hierarchical off-diagonal low-rank approximation of
  Hessians in inverse problems, with application to ice sheet model initialization}

\author{Tucker Hartland$^{1}$, Georg Stadler$^{2}$, Mauro Perego$^{3}$, Kim Liegeois$^{3}$, No\'{e}mi Petra$^{1}$}
\address{$^{1}$ Department of Applied Mathematics, University of California, Merced}
\ead{thartland@ucmerced.edu}
\address{$^{2}$ Courant Institute of Mathematical Sciences, New York University}
\address{$^{3}$ Center for Computing Research, Sandia National Laboratories}

\begin{abstract}
Obtaining lightweight and accurate approximations of Hessian applies in
inverse problems governed by partial differential equations (PDEs) is
an essential task to make both deterministic and Bayesian statistical
large-scale inverse problems computationally tractable.  The
$\mathcal{O}\left(N^3\right)$ computational complexity of dense linear
algebraic routines such as that needed for sampling from Gaussian
proposal distributions and Newton solves by direct linear methods, can be reduced to log-linear complexity by
utilizing hierarchical off-diagonal low-rank (HODLR) matrix
approximations. In this work, we show that a class of Hessians that
arise from inverse problems governed by PDEs are well approximated by
the HODLR matrix format.  In particular, we study inverse problems
governed by PDEs that model the instantaneous viscous flow of ice
sheets.  In these problems, we seek a spatially distributed basal
sliding parameter field such that the flow predicted by the ice sheet
model is consistent with ice sheet surface velocity observations.
We demonstrate the use of HODLR approximation by efficiently
generating Hessian approximations that allow
fast generation of samples from a
Gaussianized posterior proposal distribution.
Computational studies are performed
which illustrate ice sheet problem regimes for which the Gauss-Newton
data-misfit Hessian is more efficiently approximated by the HODLR
matrix format than the low-rank (LR) format. We then demonstrate that
HODLR approximations can be favorable, when compared to global
low-rank approximations, for large-scale problems by studying the
data-misfit Hessian associated to inverse problems governed by the
Stokes flow model on the Humboldt glacier and Greenland ice
sheets.
\end{abstract}



\maketitle

\section{Introduction}
\label{sec:intro}

Model-based simulation of complex physical systems plays an essential
role in understanding real world phenomena. These models are often
characterized by partial differential equations (PDEs), and are
typically subject to uncertainties stemming from unknown coefficient
fields, constitutive laws, source terms, initial and/or boundary
conditions, geometries, etc. When observation data exist, these
parameters can be estimated by solving an inverse problem governed by
the underlying model (e.g., PDE). It is well known that uncertainty is
a fundamental feature of inverse problems, therefore in addition to
inferring the parameters of interest, we need to quantify the
uncertainty associated with this inference. This uncertainty quantification
 can be done via
Bayesian inference. Solving Bayesian inverse problems governed by
complex PDEs can be extremely challenging due to high-dimensional
parameter spaces that stem from discretization of infinite-dimensional
parameter fields and the need to repeatedly solve the underlying PDEs.
To overcome these computational challenges, it is essential to
exploit problem structure, when possible. For example, the
underlying PDE solution operator is often diffusive, that observation data may be
sparse or only contain limited
information about the parameter field. These
particularities give rise to a low-rank structure in the second
derivative of the data-misfit component of the inverse problem objective (or of the negative log
likelihood), hereafter referred to as the data-misfit Hessian. In previous work~\cite{isaacpetraetal2014,
  petramartinetal2014} we exploited this low-rank structure in
the context of inverse ice sheet flow problems. However, for cases
when this ``low-rank'' is in fact large, as is the case for many inverse problems of practical interest, where the observation data are highly informative, low-rank approximation is insufficient. In this article, we exploit the local sensitivity of model predictions
to parameters, which gives rise to an off-diagonal low-rank
structure. We do so by invoking hierarchical off-diagonal low-rank
(HODLR) matrix approximations and detail how they can be used to reduce the computational
cost to solve large-scale PDE-based inverse problems.

\paragraph{Related work} 

Global low-rank approximation of Hessians in inverse problems have
been successfully utilized in~\cite{isaacpetraetal2014, SpantiniSolonenCuiEtAl15,
  flath2011, buighattasetal2013, saibaba2015}, with deterministic and randomized
methods~\cite{martinsson2011fast, buighattasetal2013}
being available to generate said approximations. However, some problems,
specifically those with highly informative observation data, are not amenable to
global low-rank approximation, and thus other structure-exploiting
strategies are needed such as those based on local
translation invariance and localized
sensitivities~\cite{algerraoetal2019, algerhartland2022,
  ZhuLiFomelEtAl16}. Here we
focus on hierarchical low-rank methods for which convenient randomized methods
are available~\cite{linlulexing2011, martinsson2016}.

Hierarchical matrices have been demonstrated in~\cite{geogaanitescustein2019, litvenkosunetal2019} to be an effective
means to approximate covariance matrices associated to large-scale
Gaussian processes. In~\cite{ambartsumyan2020}, hierarchical matrix
approximations with general hierarchical partitioning patterns are
utilized for the  construction of explicit representations of
Hessian inverses. In one of the examples studied, the authors find that
the diffusivity of the parameter-to-PDE-solution map and the informativeness of the observation data impact whether the
data-misfit Hessian is more suited for compression with hierarchical
or global low-rank formats. Here, we build on this study and focus on
a specific inverse problem arising in land ice modeling.

\paragraph{Contributions}
The main contributions of this work are as follows. (1) We motivate
the use of HODLR compression for data-misfit Hessians in
inverse problems governed by PDEs, and
present a detailed study for large-scale
ice sheet inverse problems, such as the Greenland ice sheet. (2) We
describe a strategy that leverages the fast manipulation of
HODLR matrices to efficiently generate approximate samples from a Gaussian posterior
distribution for uncertainty
quantification.
(3) We numerically
study the influence of various problem setups on the off-diagonal
low-rank structure of the data-misfit Hessian. The results show the effectiveness of the HODLR
approximation for various problem scales including for a Greenland
ice sheet inverse problem, which has a discretized parameter
dimension of $3.2\times 10^{5}$.

\section{Preliminaries}
In this section, we summarize background material regarding the
solution of discretizations of infinite-dimensional inverse problems.
We also briefly review HODLR matrices. Specifically, we define HODLR
matrices, list some of their properties and summarize the
computational complexities of computing symmetric HODLR matrix
approximations of symmetric operators that are only available through
their application on vectors. 
We refer to~\cite{hackbusch1999, hackbuschbohm2002} for a more thorough discussion
of hierarchical matrices and to~\cite{martinsson2016} for more detail on
HODLR matrices.

\subsection{Bayesian Inverse Problems}\label{sec:Bayes}
A means to account for uncertainty in parametric inference is to
employ the Bayesian approach to inverse problems~\cite{tarantola2005,
  kaipio2006, stuart2010}, which takes as input observation data $\bm{d}$, i.e., the data, prior
knowledge of the parameter and a model for the likelihood of data
conditional to $\beta$. Prior knowledge of the discretized parameter
$\bm{\beta}$ is typically determined by the expertise of domain
scientists and mathematically encoded in a probability density
function $\pi_{\text{prior}}\left(\bm{\beta}\right)$.  The likelihood
$\pi\left(\bm{d}|\bm{\beta}\right)$ involves the data uncertainty and
the mathematical model for the parameter-to-observable process.
The solution of a Bayesian inverse problem is a probability density function for the
discretized parameter $\bm{\beta}$, that is conditioned on the observation data
according to Bayes formula
\begin{equation*}
\pi_{\text{post}}\left(\bm{\beta}\right) =
\pi\left(\bm{\beta}|\bm{d}\right) \propto
\pi_{\text{prior}}\left(\bm{\beta}\right)
\pi\left(\bm{d}|\bm{\beta}\right),
\end{equation*}
which provides a formal expression for the posterior
distribution. Here, ``$\propto$'' means equal up to a normalization
constant. For a problem with Gaussian prior
$\mathcal{N}\left(\bm{\overline{\beta}},
\bm{\Gamma}_{\text{prior}}\right)$ and data noise $\bm{\eta}$
described by the zero mean Gaussian
$\mathcal{N}\left(\bm{0},\bm{\Gamma}_{\text{noise}}\right)$,
$\pi_{\text{post}}(\cdot)$ has the following form
\begin{equation}
\label{posteriorexpression}
\pi_{\text{post}}\left(\bm{\beta}\right)
\propto
\exp\left(-
\frac{1}{2}\|\bm{\mathcal{F}}(\beta)-\bm{d}\|_{\bm{\Gamma}_{\text{noise}}^{-1}}^{2}
-\frac{1}{2}\|\bm{\beta}-\bm{\overline{\beta}}
\|_{\bm{\Gamma}_{\text{prior}}^{-1}}^{2}\right),
\end{equation}
where $\bm{\mathcal{F}}$ is the parameter-to-observable map.
The notation $\|\cdot\|_{\bm{A}}$ means that the norm is weighted with the
positive-definite matrix $\bm{A}$, i.e., $\|\bm v\|_{\bm{A}}=\sqrt{\bm v^\top\bm{A}\bm v}$.
The parameter-to-PDE-solution map is typically
nonlinear, and consequently the posterior
distribution is not a Gaussian. One characteristic of the posterior
distribution is the point at which it is maximized, or equivalently the point which minimizes the negative log-posterior, the so-called maximum a posteriori (MAP) point,
\begin{equation} \label{MAPexpression}
\bm{\beta}^{\star}:=\arg\text{min}_{\bm{\beta}}
\,J(\bm{\beta}):=
\frac{1}{2}\|
\bm{\mathcal{F}}(\beta)
-\bm{d}
\|_{\bm{\Gamma}_{\text{noise}}^{-1}}^{2}
+\frac{1}{2}\|\bm{\beta}-\bm{\overline{\beta}}
\|_{\bm{\Gamma}_{\text{prior}}^{-1}}^{2}. 
\end{equation}
A means to compute the MAP point is to employ a (Gauss) Newton
method for optimization~\cite{nocedalwright2006}, which critically
relies on the availability of the (Gauss-Newton) Hessian. Since, $J$ is defined implicitly in terms of the parameter-to-observable map, which involves a PDE solution operator, we utilize the adjoint method~\cite{borzi2011, gunzburger2002, petrasachs2021} to compute it's gradient and Hessian-applies.

To fully explore posterior distributions, Markov chain Monte-Carlo
(MCMC) techniques~\cite{hastings1970, robertcasella1999} can be used.
Such techniques require a proposal distribution that ideally
approximates the posterior and is easily sampled from. One method to
generate a Gaussian proposal distribution is through the Laplace
approximation of the posterior about $\bm{\beta}_{k}$ (or around the MAP point)
\begin{equation*}
\tilde{\pi}_{\text{post}}
\left(\bm{\beta},\bm{\beta}_{k}\right)
\propto
\exp\left(-\frac{1}{2}\langle
\bm{\beta}-\bm{\mu}_{k},
\bm{H}_{k}\left(\bm{\beta}-\bm{\mu}_{k}\right)
\rangle_{\ell^{2}}\right),\\
\bm{\mu}_{k}=\bm{\beta}_{k}-\bm{H}_{k}^{-1}\bm{g}_{k},
\end{equation*}
where $\bm{g}_{k}$, $\bm{H}_{k}$ are the gradient and Hessian of the
log-posterior $J(\bm{\beta})$ at $\bm{\beta}_{k}$. Another MCMC
sampling approach is the generalized preconditioned Crank-Nicholson
(gpCN) method~\cite{rudolf2018generalization,
  pinski2015algorithms}. An attractive choice for the preconditioner
is the Hessian at the MAP point,~\cite{kim2021hippylib}.

For these and other MCMC samplers, one typically needs to apply the
inverse Hessian
$\bm{H}_{k}^{-1}$ or its square root $\bm{H}_{k}^{-1/2}$ repeatedly
and efficiently, which also motivates the study presented in this
paper. In particular, in 
Section~\ref{subsec:HODLRGaussianizedPosterior} we discuss how HODLR
approximations can be used for the fast application of the Hessian
square root.

\subsection{Symmetric HODLR Matrices}
\label{subsec:HODLRdef}

\begin{figure}[tb]
	\begin{center}  
		\begin{tikzpicture}
		\begin{scope}[yscale=0.2, xscale=0.2]
		\pgfmathtruncatemacro{\L}{1}
		\pgfmathtruncatemacro{\scalefac}{16/(2^\L)}
		\pgfmathtruncatemacro{\shiftval}{2^(\L-1)}
		\begin{scope}[yscale=-1, xscale=1]
		\begin{scope}[yshift=-\shiftval, xshift=-\shiftval]
		\begin{scope}[yscale=\scalefac, xscale=\scalefac]
		\pgfmathtruncatemacro{\msize}{2^\L}
		\pgfmathtruncatemacro{\colstep}{80}
		\foreach \l in {1,...,\L}
		{
			\pgfmathtruncatemacro{\delI}{\msize*2^(-\l)}
			\pgfmathtruncatemacro{\DelI}{2*\delI}
			\pgfmathtruncatemacro{\colorl}{75}
			\pgfmathtruncatemacro{\maxi}{2^(\l-1)}
			\foreach \i in {1,...,\maxi}
			{
				\pgfmathtruncatemacro{\a}{(\i-1)*\DelI}
				\pgfmathtruncatemacro{\b}{(\i-1)*\DelI+\delI}
				\pgfmathtruncatemacro{\c}{\i*\DelI}
				\filldraw[fill=teal] (\a, \b) rectangle (\b, \c);
				\filldraw[fill=teal] (\b, \a) rectangle (\c, \b);
			}
		}
		\pgfmathtruncatemacro{\delI}{\msize*2^(-\L)}
		\pgfmathtruncatemacro{\DelI}{2*\delI}
		\pgfmathtruncatemacro{\colorl}{20+(\L-1)*\colstep}
		\pgfmathtruncatemacro{\maxi}{2^(\L-1)}
		\foreach \i in {1,...,\maxi}
		{
			\pgfmathtruncatemacro{\a}{(\i-1)*\DelI}
			\pgfmathtruncatemacro{\b}{(\i-1)*\DelI+\delI}
			\pgfmathtruncatemacro{\c}{\i*\DelI}
			\filldraw[fill=violet] (\a, \a) rectangle (\b, \b);
			\filldraw[fill=violet] (\b, \b) rectangle (\c, \c);
		}
		\end{scope}
		\end{scope}
		\end{scope}
		\end{scope}
		\end{tikzpicture}
		\qquad
		\begin{tikzpicture}
		\begin{scope}[yscale=0.2, xscale=0.2]
		\pgfmathtruncatemacro{\L}{2}
		\pgfmathtruncatemacro{\scalefac}{16/(2^\L)}
		\pgfmathtruncatemacro{\shiftval}{2^(\L-1)}
		\begin{scope}[yscale=-1, xscale=1]
		\begin{scope}[yshift=-\shiftval, xshift=-\shiftval]
		\begin{scope}[yscale=\scalefac, xscale=\scalefac]
		\pgfmathtruncatemacro{\msize}{2^\L}
		\pgfmathtruncatemacro{\colstep}{40/(\L-1)}
		\foreach \l in {1,...,\L}
		{
			\pgfmathtruncatemacro{\delI}{\msize*2^(-\l)}
			\pgfmathtruncatemacro{\DelI}{2*\delI}
			\pgfmathtruncatemacro{\colorl}{75}
			\pgfmathtruncatemacro{\maxi}{2^(\l-1)}
			\foreach \i in {1,...,\maxi}
			{
				\pgfmathtruncatemacro{\a}{(\i-1)*\DelI}
				\pgfmathtruncatemacro{\b}{(\i-1)*\DelI+\delI}
				\pgfmathtruncatemacro{\c}{\i*\DelI}
				\ifnum\l=1
				\filldraw[fill=teal] (\a, \b) rectangle (\b, \c);
				\filldraw[fill=teal] (\b, \a) rectangle (\c, \b);
				\else
				\ifnum\l=2 
				\filldraw[fill=olive] (\a, \b) rectangle (\b, \c);
				\filldraw[fill=olive] (\b, \a) rectangle (\c, \b);
				\else
				\filldraw[fill=magenta] (\a, \b) rectangle (\b, \c);
				\filldraw[fill=magenta] (\b, \a) rectangle (\c, \b);
				\fi
				\fi 
			}
		}
		\pgfmathtruncatemacro{\delI}{\msize*2^(-\L)}
		\pgfmathtruncatemacro{\DelI}{2*\delI}
		\pgfmathtruncatemacro{\colorl}{20+(\L-1)*\colstep}
		\pgfmathtruncatemacro{\maxi}{2^(\L-1)}
		\foreach \i in {1,...,\maxi}
		{
			\pgfmathtruncatemacro{\a}{(\i-1)*\DelI}
			\pgfmathtruncatemacro{\b}{(\i-1)*\DelI+\delI}
			\pgfmathtruncatemacro{\c}{\i*\DelI}
			\filldraw[fill=violet] (\a, \a) rectangle (\b, \b);
			\filldraw[fill=violet] (\b, \b) rectangle (\c, \c);
		}
		\end{scope}
		\end{scope}
		\end{scope}
		\end{scope}
		\end{tikzpicture} 
		\qquad
		\begin{tikzpicture}
		\begin{scope}[yscale=0.2, xscale=0.2]
		\pgfmathtruncatemacro{\L}{3}
		\pgfmathtruncatemacro{\scalefac}{16/(2^\L)}
		\pgfmathtruncatemacro{\shiftval}{2^(\L-1)}
		\begin{scope}[yscale=-1, xscale=1]
		\begin{scope}[yshift=-\shiftval, xshift=-\shiftval]
		\begin{scope}[yscale=\scalefac, xscale=\scalefac]
		\pgfmathtruncatemacro{\msize}{2^\L}
		\pgfmathtruncatemacro{\colstep}{40/(\L-1)}
		\foreach \l in {1,...,\L}
		{
			\pgfmathtruncatemacro{\delI}{\msize*2^(-\l)}
			\pgfmathtruncatemacro{\DelI}{2*\delI}
			\pgfmathtruncatemacro{\colorl}{75}
			\pgfmathtruncatemacro{\maxi}{2^(\l-1)}
			\foreach \i in {1,...,\maxi}
			{
				\pgfmathtruncatemacro{\a}{(\i-1)*\DelI}
				\pgfmathtruncatemacro{\b}{(\i-1)*\DelI+\delI}
				\pgfmathtruncatemacro{\c}{\i*\DelI}
				\ifnum\l=1
				\filldraw[fill=teal] (\a, \b) rectangle (\b, \c);
				\filldraw[fill=teal] (\b, \a) rectangle (\c, \b);
				\else
				\ifnum\l=2 
				\filldraw[fill=olive] (\a, \b) rectangle (\b, \c);
				\filldraw[fill=olive] (\b, \a) rectangle (\c, \b);
				\else
				\filldraw[fill=gray] (\a, \b) rectangle (\b, \c);
				\filldraw[fill=gray] (\b, \a) rectangle (\c, \b);
				\fi
				\fi 
			}
		}
		\pgfmathtruncatemacro{\delI}{\msize*2^(-\L)}
		\pgfmathtruncatemacro{\DelI}{2*\delI}
		\pgfmathtruncatemacro{\colorl}{20+(\L-1)*\colstep}
		\pgfmathtruncatemacro{\maxi}{2^(\L-1)}
		\foreach \i in {1,...,\maxi}
		{
			\pgfmathtruncatemacro{\a}{(\i-1)*\DelI}
			\pgfmathtruncatemacro{\b}{(\i-1)*\DelI+\delI}
			\pgfmathtruncatemacro{\c}{\i*\DelI}
			\filldraw[fill=violet] (\a, \a) rectangle (\b, \b);
			\filldraw[fill=violet] (\b, \b) rectangle (\c, \c);
		}
		\end{scope}
		\end{scope}
		\end{scope}
		\end{scope}
		\end{tikzpicture}
	\end{center}
	\caption{Rank-structure of a
		matrix $\bm{A}$ with hierarchical depths $L=1$ (left), $L=2$ (middle)
		and $L=3$ (right). Off-diagonal blocks are assumed to be low-rank.}  
	\label{fig:hmatrixpartitioningstructure} 
\end{figure}
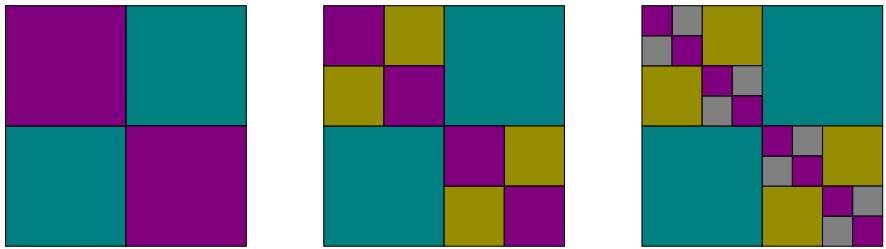

%
A HODLR matrix $\bm{A}\in\mathbb{R}^{N \times N}$, is a matrix
equipped with a depth $L\in\mathbb{N}$, hierarchical partitionings of
the index set $\mathcal{I}=\lbrace 1,2,\dots,N\rbrace$ into continguous subsets and low-rank
off-diagonal blocks defined by the partition, which is described in greater detail in e.g.~\cite{martinsson2016}. The block rank-structure of a HODLR matrix for various hierarchical depths is illustrated in Figure~\ref{fig:hmatrixpartitioningstructure}. An HODLR matrix must satisfy two additional properties.
\begin{enumerate}
	\item The depth of the hierarchical
	partitioning scales with the logarithm
	of the size of the matrix, i.e.,
	\begin{equation*}
		L=\mathcal{O}\left(\log\,N\right).
	\end{equation*}
	\item The maximum rank of each hierarchical level $\ell$ off-diagonal block, $r_{\ell}$, is bounded above by a number $r$ that is independent of the problem size $N$, for each level $\ell$ 
	\begin{equation*}
		\max_{1 \leq \ell \leq L}r_{\ell}\leq r=\mathcal{O}\left(1\right).
	\end{equation*}
\end{enumerate}
Such matrices are referred to as data-sparse since
the low-rank blocks allow for them to be represented computationally with less than $\mathcal{O}\left(N^{2}\right)$ floating point
numbers. In particular, the storage of an HODLR matrix is 
$\mathcal{O}\left(N\,\log\,N\right)$,
$\mathcal{O}(N\,\log\,N)$ flops are needed
to compute a HODLR matrix-vector product~\cite{martinsson2011fast},
and $\mathcal{O}(N\,\log^{2}\,N)$ flops are required for
direct methods to compute an inverse HODLR matrix-vector product~\cite{ambikasarandarve2013}, as well as square root and inverse square root matrix-vector products~\cite{ambikasaranoneil2014}.

\paragraph{Compression}

We aim to generate HODLR approximations
of data-misfit Hessians in inverse problems. For large-scale problems, the data-misfit Hessian is available only as a matrix-free operator. In order to construct HODLR approximations of symmetric matrix-free operators,
we employ previously developed randomized linear algebraic routines
which only require the matrix-free action on a limited number of random vectors
with specified null entries, referred to
as \textit{structured} random vectors. 
The Hessian action on these structured random vectors
is used to sample row and column spaces of off-diagonal Hessian submatrices
and allow for randomized approximate truncated singular value decompositions of the aforementioned off-diagonal submatrices. More details can be found in the appendix, see Algorithm~\ref{alg:HODLRcompression}.

For the results that we present in Section~\ref{sec:2DResults}
a rank-adaptive symmetric matrix-free~\cite{halkomartinsson2011, xixiachan2014}, hierarchical compression algorithm
is utilized, that is based on~\cite{martinsson2016}.
A similar algorithm is presented in~\cite{keyesturkiyyah2019},
wherein the hierarchical partitioning
is more general and the low-rank blocks have nested bases.
The rank-adaptivity provides a high probability means of resolving
the off-diagonal blocks to a desired level
of accuracy. By utilizing available matrix-vector product information and the Rayleigh quotient, a rank adaptive relative tolerance algorithm is made possible.

\paragraph{Computational Cost of Generating HODLR Approximations}

The number of matrix-vector products $\zeta$,
needed to compress a symmetric
matrix using $d$ oversampling vectors,
into a level $L$ HODLR matrix with off-diagonal ranks
$\lbrace r_{\ell}\rbrace_{\ell=1}^{L}$ is given by
\begin{equation}
\zeta=2\left(\langle r\rangle +d\right)L+N/2^{L},
\text{ where }\langle r\rangle:=\frac{1}{L}\sum_{\ell=1}^{L}r_{\ell}. 
\label{eq:HODLRcost}
\end{equation}
Equation~\ref{eq:HODLRcost} can be understood from Algorithm~\ref{alg:HODLRcompression} in Appendix~\ref{subsec:randomizedcompressionalgorithms}, 
as for each level $\ell$ one needs to compute $r_{\ell}+d$ Hessian vector products, in order
to compute $\bm{Y}$ (line~$7$ of Algorithm~\ref{alg:HODLRcompression}) and $r_{\ell}+d$ Hessian vector products to compute
$\bm{Z}$ (line~$14$ of Algorithm~\ref{alg:HODLRcompression}). The remaining $N/2^{L}$ Hessian vector products arise from the need to determine the diagonal subblocks, which is detailed in~\cite{martinsson2011fast}. We note that with an adaptive procedure to determine an
approximate basis $\bm{Q}$, such as that in~\cite{xixiachan2014}, for a
block matrix column space, the cost is reduced to
$\zeta_{\text{adaptive}}=2\left(\langle r\rangle+d/2\right)L+N/2^{L}$
but with the additional computational burden of extra orthogonalization routine calls. We note that $\zeta=\mathcal{O}(\log\,N)$
matrix-vector products are needed to generate an HODLR approximation of a matrix with HODLR structure. For sufficiently large problems HODLR compression is not expected to be more computationally efficient than global low-rank (LR) compression, as $\zeta^{\text{LR}}=r+d$, the number of matrix-vector products to generate a rank $r$ compression by the single-pass algorithm~\cite{martinsson2016} with $d$ oversampling vectors is independent of the problems size. However, for problems of substantial size, we observe that the HODLR format does offer computational savings (see Section~\ref{sec:HumboldtGreenland}).

\section{HODLR matrices in inverse problems governed by PDEs}
\label{sec:HODLRapplication}
Here, we illustrate why data-misfit Hessians in inverse problems
governed by PDEs may contain numerically low-rank off-diagonal
blocks, describe how one can permute parameters to expose
this HODLR structure, and show how HODLR approximations can be
leveraged to draw samples from Gaussian approximations of Bayesian
posterior distributions.

\subsection{Motivation}
\label{subsec:motivation}

\begin{figure}[tb]
	\begin{center} 
		\begin{tikzpicture}
		\draw[red,dashed] (7,.5) -- (3,.5) -- (1,0);
		\draw[dashed] (3,.5) -- (3,2);
		\draw[blue,dashed] (3,2) -- (3,2.5);
		\fill[red!90,nearly transparent] (5,0) -- (7,.5) -- (3,.5) -- (1,0) -- cycle;
		\draw [fill=red!40!white,opacity=1] (3,2.25) -- (3.5,2.25) -- (3.28,.25) -- (3.22,0.25) -- cycle;
		\draw (3.25,0.25) -- (2.25, -0.5);
		\draw [thick] (2.25, -0.5) node[below]{perturbation $\psi_{i}$};
		\draw [fill=red] (3.25,2.25) circle (.25cm and 0.07cm);
		\draw [fill=red] (3.25,0.25) circle (.025cm and 0.007cm);
		\draw (3.25,2.25) -- (2,3);
		\draw [thick] (2,3) node[above]{sensitivity cone, $\frac{\delta u}{\delta \beta}(\beta)(\psi_{i})$};
		\draw (1,0) -- (5,0) -- (5,2) -- (1,2) -- (1,0);
		\draw  (5,2) -- (7,2.5) -- (3,2.5) -- (1,2);
		\draw  (7,2.5) -- (7,.5) -- (5,0);
		\draw  (5,2.25) -- (6,3) ;
		\draw [thick] (6,3) node[above]{$\Gamma_{2}$};
		\draw  (4.75,.125) -- (5.75,-.5) ;
		\draw [thick] (5.75,-.5) node[below]{$\Gamma_{1}$};
		\fill[blue!90,nearly transparent] (5,2) -- (7,2.5) -- (3,2.5) -- (1,2) -- cycle;
		\end{tikzpicture}
	\end{center}
	\caption{Sketch illustrating a case where the influence of
          changes in
          the parameter $\beta$ on the PDE solution $u$ in
          $\Gamma_{2}$ is focused in a small area. To
          illustrate this, we show a sensitivity cone,
          i.e., the PDE
          solution $u$ is predominantly impacted in a
          cone about the support of the localized
          parameter perturbation.}
	\label{fig:sensitivitycone} 
\end{figure}
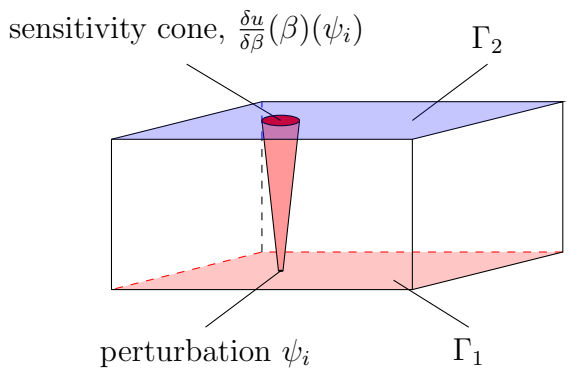

Consider the following data-misfit cost functional
\begin{equation*}
J_{\text{misfit}}\left(\beta\right):=\frac{1}{2}\|\bm{\mathcal{F}}(\beta)-\bm{d}\|
_{\bm{\Gamma}_{\text{noise}}^{-1}}^{2},\quad \text{with}\quad \mathcal{F}(\beta)=\bm{\mathcal{B}}u,
\end{equation*} 
where $\bm{\mathcal{B}}$ linearly maps the PDE solution $u=u(\beta)$,
for the spatially-distributed parameter field $\beta$, to the model
predictions associated to the data $\bm{d}$. Moreover,
$\bm{\Gamma}_{\text{noise}}$ is the covariance matrix describing the Gaussian
noise of the observational data. For illustration purposes, we assume
that the parameter function $\beta$ is defined on a region $\Gamma_{1}$ and the data $\bm{d}$ is
observed on a region $\Gamma_{2}$, which may or may not be distinct.
These quantities are related through the solution of the governing PDE
and the measurement operator $\bm{\mathcal B}$. The
characteristics of this relation depends on properties of the
governing PDE. In
the following, we assume that a spatially (or temporally) localized
perturbation in the $\beta$ field leads to a predominantly localized
effect in the PDE solution $u$, and thus the model predictions
$\bm{\mathcal{B}}u$.  This property is illustrated in Figure~\ref{fig:sensitivitycone}, where we use a sensitivity cone to
illustrate the influence of a local perturbation in
$\beta$, defined over $\Gamma_1$, on the PDE solution $u$ in
$\Gamma_2$. It is well known that for an elliptic PDE, local
perturbations influence the solution globally, but depending on the
geometry of the domain and the equation, this global effect may
rapidly decay outside a subset of $\Gamma_2$ that captures the main
effects of the perturbation. For instance, in a problem as in Figure~\ref{fig:sensitivitycone}, the influence of perturbations in $\beta$
on $u$ is likely to become more localized when the distance between
$\Gamma_1$ and $\Gamma_2$ decreases.

%

We next discuss the relationship between properties of the PDE as
discussed above and off-diagonal blocks in the Hessian matrix (or its
Gauss-Newton variant). The data-misfit Hessian, i.e., the Hessian of
the data-misfit part of the cost functional, can be derived using the
adjoint method~\cite{borzi2011, gunzburger2002, petrasachs2021}. However, we find
that the HODLR structure of the data-misfit Hessian is most easily
seen by studying a formal expression of it in terms of the first and
second order sensitivities $\delta u/\delta\beta$,
$\delta^{2}u/\delta\beta^{2}$ 
%
\begin{eqnarray*}
\frac{\delta^{2}}{\delta\beta^{2}}J_{\text{misfit}}\left(\beta\right)\left(\beta_{1},\beta_{2}\right)
&=&
\left(\bm{\mathcal{B}}u-\bm{d}\right)^{\top}
\bm{\Gamma}_{\text{noise}}^{-1}
\left(
\bm{\mathcal{B}}
\frac{\delta^{2}u}{\delta \beta^{2}}\left(\beta\right)\left(\beta_{1},\beta_{2}\right)\right)\,
+ \\
&\phantom{e}&
\left(
\bm{\mathcal{B}}
\frac{\delta u}{\delta \beta}\left(\beta\right)\left(\beta_{1}\right) 
\right)^{\top}\,
\bm{\Gamma}_{\text{noise}}^{-1}
\left(
\bm{\mathcal{B}}
\frac{\delta u}{\delta \beta}\left(\beta\right)\left(\beta_{2}\right) 
\right),  
\end{eqnarray*}
 where $\delta u/\delta \beta\left(\beta\right)\left(\beta_{1}\right)$
 is the first variation~\cite{gelfand1963} of $u$ with respect to
 $\beta$ in direction $\beta_{1}$, and $\delta^{2}u/\delta
 \beta^{2}\left(\beta\right)\left(\beta_{1},\beta_{2}\right)$ is the
 second variation of $u$ with respect to $\beta$ in directions
 $\beta_{1},\beta_{2}$, that is,
 \begin{eqnarray*}
\frac{\delta u}{\delta\beta}\left(\beta\right)\left(\beta_{1}\right):=
\left[\frac{\mathrm{d}}{\mathrm{d}\epsilon}
u\left(\beta+\epsilon\beta_{1}\right)\right]_{\epsilon=0}, \\
\frac{\delta^{2}u}{\delta\beta^{2}}\left(\beta\right)\left(\beta_{1},\beta_{2}\right):=
\left[\frac{\mathrm{d}}{\mathrm{d}\epsilon}
\frac{\delta u}{\delta\beta}\left(\beta+\epsilon\beta_{2}\right)\left(\beta_{1}\right)\right]_{\epsilon=0}.
\end{eqnarray*}
Upon discretizing $\beta$ with finite elements we obtain the following
formal expression for the $(i,j)$-entry of the data-misfit Hessian
$\bm{H}_{\text{misfit}}$ and of the Gauss-Newton
data-misfit Hessian $\bm{H}_{\text{misfit}}^{\text{GN}}$
\begin{eqnarray}
\label{eq:Hmisfitelements}
\left(\bm{H}_{\text{misfit}}\right)_{i,j}&=
\frac{\delta^{2}}{\delta \beta^{2}}\Big(J_{\text{misfit}}\left(\beta\right)\Big) 
\left(\psi_{i},\psi_{j}\right),\\
\left(\bm{H}_{\text{misfit}}^{\text{GN}}\right)_{i,j}&=
\left(
\bm{\mathcal{B}}
\frac{\delta u}{\delta \beta}\left(\beta\right)\left(\psi_{i}\right) 
\right)^{\top} \,
\bm{\Gamma}_{\text{noise}}^{-1}
\left(
\bm{\mathcal{B}}
\frac{\delta u}{\delta \beta}\left(\beta\right)\left(\psi_{j}\right) 
\right),
\end{eqnarray}
where $\lbrace \psi_{j} \rbrace_{j=1}^{N}$ is a basis for the nodal
finite-element space, which is used to approximate~$\beta$.
%
%

When sensitivities are predominantly local as discussed above and when the support of two finite
element basis functions $\psi_{i},\psi_{j}$ are well separated, the terms
\begin{equation*}
	\left(
	 \bm{\mathcal{B}}
	\frac{\delta u}{\delta \beta}\left(\beta\right)\left(\psi_{i}\right)
	\right)^{\top}\,
	\bm{\Gamma}_{\text{noise}}^{-1}
	\Big(\bm{\mathcal{B}}\frac{\delta u}{\delta \beta}\left(\beta\right)\left(\psi_{j}\right) \Big)
	\quad\text{and}\quad
	\bm{\mathcal{B}}\Big(\frac{\delta^{2}u}{\delta \beta^{2}}\left(\beta\right)\left(\psi_{i},\psi_{j}\right)\Big),
\end{equation*}
are rather small (assuming diagonally dominant noise covariance
matrices). This is, e.g., due to $\bm{\mathcal{B}}\delta
u/\delta\beta(\beta)(\psi_{i})$ having small values when
$\bm{\mathcal{B}}\delta u/\delta \beta(\beta)(\psi_{j})$ is
large. Now, let $\mathcal{I},\mathcal{J}$ be disjoint index
subsets of $\lbrace 1,2,\dots,N\rbrace$, then the entries in the matrix
block $\lbrace \left(\bm{H}_{\text{misfit}}\right)_{i\in \mathcal
  I,j\in \mathcal J}\rbrace$ of the data-misfit Hessian are relatively
small whenever
$\cup_{i\in\mathcal{I}}\text{supp}\left(\psi_{i}\right)$ and
$\cup_{j\in\mathcal{J}}\text{supp}\left(\psi_{j}\right)$ are well
separated. Such Hessian blocks are well suited for
approximation by low-rank matrices. When the degrees of freedom
corresponding to the finite element basis functions $\psi_i$
are ordered such that
$\mathcal{I},\mathcal{J}$ are contiguous,
$\left(\bm{H}_{\text{misfit}}\right)_{\mathcal{I},\mathcal{J}}$ is an
off-diagonal subblock of $\bm{H}_{\text{misfit}}$ and
$\bm{H}_{\text{misfit}}$ tends to have HODLR structure as defined in
Section~\ref{subsec:HODLRdef}. The Gauss-Newton data-misfit Hessian
may have HODLR structure for the same reasons. In both cases, the
order of the basis functions and thus the degrees of freedom influence
this structure. Ideally, one wants an order that maintains locality,
i.e., consecutive indices correspond to basis functions that are close
to each other, and as a consequence, basis function with significantly
different indices are far from each other such that the corresponding
off-diagonal blocks have small entries and can be well approximated
using a low-rank matrix approximation.  We defer to Section~\ref{subsec:dofordering} for a discussion of methods and numerical
experiments regarding the order of the degrees of freedom.

\subsection{Application of HODLR structure for fast sampling of
  Gaussian posterior approximations}
\label{subsec:HODLRGaussianizedPosterior}

In \cite{petramartinetal2014},
the following expressions of the 
Gaussianized posterior covariance are provided,
\begin{eqnarray*}
\bm{\Gamma}_{\text{post}}&=&
\left(\bm{H}_{\text{misfit}}+
\bm{\Gamma}_{\text{prior}}^{-1}\right)^{-1} =
\bm{\Gamma}_{\text{prior}}^{1/2}
\left(
\bm{H}_{\text{misfit}}^{\prime}
+\bm{I}
\right)^{-1}
\bm{\Gamma}_{\text{prior}}^{\top/2}, \\
\bm{H}_{\text{misfit}}^{\prime}&:=&
\bm{\Gamma}_{\text{prior}}^{\top/2}
\bm{H}_{\text{misfit}}
\bm{\Gamma}_{\text{prior}}^{1/2}, \\			
\bm{\Gamma}_{\text{post}}^{1/2}
&=&\bm{\Gamma}_{\text{prior}}^{1/2}
\left(
\bm{H}_{\text{misfit}}^{\prime}
+\bm{I}
\right)^{-1/2},
\end{eqnarray*}
where the matrix square-root $\bm{A}^{1/2}$ is such that
$\bm{A}=\bm{A}^{1/2}\left(\bm{A}^{1/2}\right)^{\top}$.  For Bayesian
inverse problems with a parameter field that is distributed spatially
over a bounded subset of $\mathbb{R}^{m}$, $m=2,3$, a reasonable
choice is to use the square of an inverse elliptic PDE operator for
the prior covariance~\cite{stuart2010}, which permits a means of
obtaining a symmetric square root of
$\bm{\Gamma}_{\text{prior}}$. In previous works such as~\cite{isaacpetraetal2014, SpantiniSolonenCuiEtAl15, flath2011, buighattasetal2013, saibaba2015},
the prior-preconditioned data-misfit Hessian
$\bm{H}_{\text{misfit}}^{\prime}$, was approximated by global low-rank
compression. This strategy provides an efficient means of approximating the posterior covariance matrix in inverse problems with data sets that contain sufficiently small amounts of information.
Here we propose to exploit HODLR problem structure and generate approximate posterior covariance matrices by HODLR approximations of the prior-preconditioned data-misfit $\bm{\tilde{H}}_{\text{misfit}}^{\prime}$,
see Appendix~\ref{subsec:posteriorerroranalysis} for an analysis on
how such an approximation impacts the accuracy of the approximate
posterior covariance
\begin{equation*} 
\bm{\tilde{\Gamma}}_{\text{post}}=\bm{\Gamma}_{\text{prior}}^{1/2}
\left(
\bm{\tilde{H}}_{\text{misfit}}^{\prime}
+\bm{I}
\right)^{-1}
\bm{\Gamma}_{\text{prior}}^{\top/2}.
\end{equation*}
A symmetric square-root factorization
of $\bm{\tilde{H}}_{\text{misfit}}^{\prime}+\bm{I}$ is then generated
with $\mathcal{O}\left(N\,\log^{2}\,N\right)$
flops~\cite{ambikasaranoneil2014}. 
The symmetric factorization allows for
a $\mathcal{O}\left(N\,\log\,N\right)$ 
means of applying both the square root
and inverse square root.

\section{Bayesian inverse ice sheet problems}

The simulation of the dynamics of ice sheets (e.g., the Greenland or
Antarctic ice sheets) is an important component of coupled climate
simulations.  Such simulations require estimation of a
present state of the ice that is consistent with available
observations, a process sometimes referred to as model
initialization. This estimation problem can be formulated either as a
deterministic inverse problem (i.e., as nonlinear least squares
optimization governed by PDEs) or as a Bayesian inverse problem (i.e., as a statistical
problem which aims to characterize a distribution of states). The
latter approach, while more expensive, provides uncertainty estimates
in addition to determining a best parameter fit.

Ice sheet dynamics~\cite{cuffeypatterson2010} is typically governed by nonlinear Stokes equations
or simplifications thereof, such as the first-order equations 
(see e.g.,~\cite{dukowicz2010}). Generally, the most uncertain component in
ice sheet simulations is the basal boundary condition, i.e., how the
ice sheet interacts with the rock, sand, water or a mix thereof at its
base. Estimating an ice sheet's effective boundary condition from
velocity observations on the top surface, the ice sheet's geometry and
a model for its dynamics is thus an important problem that can
mathematically formulated as an inverse problem~\cite{isaacpetraetal2014, 
	larouretal2012,
	morlighem2010,
	peregopricestadler2014, 
	petrazhustadleretal2012}.

We summarize the formulation of this inverse problem next. As common
in the literature, we use \emph{a snapshot} optimization approach, where all the data are assumed to be collected over a short period of time during which changes in the ice geometry are negligible. We denote
the bounded domain covered by ice by $\Omega\subset\mathbb{R}^{m}$, $m\in\lbrace 2,3\rbrace$,
and the basal, lateral and top parts of the domain boundary
$\partial\Omega$ by $\Gamma_{b}$, $\Gamma_{l}$, and $\Gamma_{t}$, as illustrated in Figure~\ref{fig:schematic}.

The governing equations are nonlinear incompressible Stokes
equations whose solution is the ice flow velocity $\bm{u}:\Omega\to\mathbb{R}^m$ and the pressure
$p:\Omega\to\mathbb{R}$ given as follows:
\begin{eqnarray}
	-\nabla\cdot\bm{\sigma}_{\bm{u}}
	=\rho\bm{g}\,\,\,\text{ in }\Omega, \label{Stokeseqn:1}\\
	\phantom{-}\nabla\cdot\bm{u}
	=0\,\,\,\,\,\,\,\,\,\,\text{ in }\Omega, \label{Stokeseqn:2}\\
	\phantom{-}\bm{\sigma}_{\bm{u}}\bm{n}
	=\bm{0}\,\,\,\,\,\,\,\,\,\,\text{ on }\Gamma_{t}, \label{Stokeseqn:3}\\
	\phantom{-}\bm{u}\cdot\bm{n} =0 \text{ and } \bm{T}\left(\bm{\sigma}_{\bm{u}}\bm{n}+\exp\left(\beta\right)\bm{u}\right)
	=\bm{0}\,\,\,\,\,\,\text{ on }\Gamma_{b},\label{Stokeseqn:4}
\end{eqnarray}
%
%
along with additional lateral boundary conditions. Here, $\beta$ is a
basal sliding parameter field, $\rho\bm{g}$ the body force density,
where $\rho$ is the mass density of the ice and $\bm{g}$ the
acceleration due to gravity. Equation~\ref{Stokeseqn:1}
describes the conservation of momentum,~\ref{Stokeseqn:2} the
conservation of mass, and~\ref{Stokeseqn:3} are stress-free boundary
conditions for the top surface (the ice-air interface). In normal
direction, Equation~\ref{Stokeseqn:4} states a non-penetration condition,
i.e., the ice cannot flow into the rock/sand layer which supports it (here
$\bm{n}$ denotes the outward unit normal to the boundary $\partial\Omega$ and $\bm{T}$ the tangential operator, $\bm{Tv} = \bm{v}-\bm{n}(\bm{n}^{\top}\bm{v})$).  In tangential
direction, Equation \ref{Stokeseqn:4} specifies a tangential sliding
condition that relates the fraction of tangential sliding and
tangential stress through the (logarithmic) basal sliding field
$\beta=\beta(x)$, $x\in \Gamma_b$.  We employ Glen's flow law~\cite{glen1955}, a constitutive law for ice that relates the stress
tensor $\bm{\sigma}_{\bm{u}}$ and the strain rate tensor
$\bm{\dot{\varepsilon}}_{\bm{u}}= \frac{1}{2}\left(
\bm{\nabla}\bm{u}+\bm{\nabla}\bm{u}^{\top} \right)$,
\begin{equation}
\bm{\sigma}_{\bm{u}}= 2\eta\left(\bm{u}\right)
\bm{\dot{\varepsilon}}_{\bm{u}} -\bm{I}p, \text{ with } \eta\left(\bm{u}\right) =
\frac{1}{2}A^{-1/n}\bm{\dot{\varepsilon}}_{\text{II}}^{\frac{1-n}{2n}},
\end{equation}
where $\eta$ is the effective viscosity, $\bm{I}$ is the unit matrix,
$\bm{\dot{\varepsilon}}_{\text{II}} =
\text{tr}\left(\bm{\dot{\varepsilon}}_{\bm{u}}^{2}\right)$ is the
second invariant of the strain rate tensor, $A$ is a flow rate factor,
and $n$ is Glen's exponent. Ice is typically modeled using $n\approx
3$, which corresponds to a shear-thinning constitutive relation, here we use $n=3$.

As discussed above, the parameter containing the largest uncertainty
is the (logarithmic) basal sliding field $\beta=\beta(x)$. Thus, it is
usually the parameter inferred from (typically, satellite)
observation data $\bm{d}$, here in the form of surface velocity measurements. Using an
appropriate point observation operator $\bm{ \mathcal B}$ that extracts point
data from the solution $\bm{u}$ of the governing equations~\ref{Stokeseqn:1}-\ref{Stokeseqn:4}, and assuming additive observation errors $\bm{\eta}$,
the relationship between model and data is now of the typical
form
\begin{equation}
\bm d = \bm{\mathcal B}\bm u + \bm \eta.
\end{equation}
Assuming that the observation errors $\bm \eta$ and the prior
for the parameter field $\beta$ follow Gaussian distributions, we are in the
framework of Bayesian inverse problems summarized in
Section~\ref{sec:Bayes}.

\section{Example I: Two-dimensional ISMIP-HOM benchmark}
\label{sec:2DResults}
We first study the prospects of compressing the Gauss-Newton data-misfit Hessian in a
problem inspired by the ISMIP-HOM collection of ice sheet simulation
benchmark problems~\cite{pattynetal2008}. This problem set
was used to explore inverse ice sheet problems 
in~\cite{peregopricestadler2014, petrazhustadleretal2012}.
After a short description of the problem setup, we present results such as
the MAP point estimate $\bm{\beta}^{\star}$ and approximate Gaussianized posterior samples using an HODLR compression of the posterior covariance. Then, we study the impact that various
problem features have on the suitability of the Gauss-Newton data-misfit Hessian for compression
to the HODLR and global low-rank formats.

\subsection{Problem setup}
\begin{figure}[tb]
	\centering
	\begin{tikzpicture}
	\draw[->,thick] (-1,-1/6*-1) -- (7,-1/6*7) node[pos=0.97, above] {${ x}$} node [pos=0.33, below] {$\Gamma_{b}$} ;
	\draw[->,thick] (0,0) -- (1/6*3,3) node[pos=0.93, left]  {${z}$};
	\shade[bottom color=lightblue,nearly transparent]  (0,0) -- (1/6*2,2)  -- (1/6*2+6,2-1/6*6) -- (6,-1/6*6) -- (0,0);
	\shade[bottom color=gray,shading angle=60,nearly transparent]  (-1,-1/6*-1) -- (7,-1/6*7) -- (-1,-1/6*7)-- (-1,-1/6*-1);
	\draw[thick] (0,0) -- (1/6*2,2)  -- (1/6*2+6,2-1/6*6) -- (6,-1/6*6) -- (0,0);
	\draw[dashed,thick] (6,-1/6*6) -- (3,-1/6*6);
	\path[->,thick] (4,-4/6) edge [bend right] (4,-1);
	\node at (3.7,-5/6) {$\theta$};
	\draw [<->,thick] (-.5-.1,.5*1/6) -- (1/6*2-.5-.1,2+.5*1/6) node [midway, left] {$H$} node [midway, right] {$\Gamma_{l}$};
	\node[right] at (6+1/6,0) {$\Gamma_{l}$};
	\draw [<->,thick] (1/6*2+6+.5/6,2-1/6*6+6*.5/6+.1) --  (1/6*2+.5/6,2+6*.5/6+.1) node [midway, above] {$W$} node [midway, below] {$\Gamma_{t}$};
	\foreach \x in {1,2,...,20}{
		\pgfmathsetmacro{\y}{1/2+1/2*rand}
		\draw [thick,blue]({1/6*2*\y+(1-\y)*(1/6*2+6)},{2*\y+(2-1/6*6)*(1-\y)}) circle [radius=0.05];
	}
	\end{tikzpicture}
	\caption{Schematic of two-dimensional slab of ice used for
		Example I in
		Section~\ref{sec:2DResults}. The blue circles show
		representative (random) measurement locations. The angle $\theta$ is the slope of the ice
		slab.}\label{fig:schematic}
\end{figure}
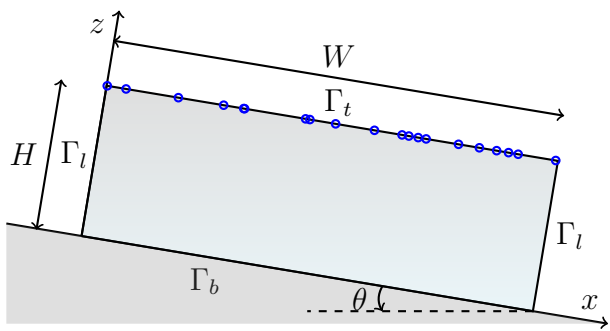

This problem setup consists of a rectangular piece of ice on a
slope, as sketched in Figure~\ref{fig:schematic}. This simple example allows
us to study the influence of the domain aspect ratio, the number of observations and the level of mesh refinement on the properties of the Gauss-Newton data-misfit Hessian matrix.
The domain has a width of
$W=10^{4} \left[\text{m}\right]$ and a height of
$H=10^{2}\left[\text{m}\right]$. Periodic boundary conditions are
employed along the lateral boundaries such that the setup models an
infinite slab of ice on a slope. The governing equations and
other boundary conditions are as discussed in 
Equations~\ref{Stokeseqn:1}-\ref{Stokeseqn:4}.

The Stokes equations are discretized using Taylor-Hood finite elements
on a mesh of~$256\times 10$ rectangles, each subdivided into two
triangles, for the domain length $\left[0,W\right)$ and height
$\left[0,H\right]$.  To compute a MAP estimate, we generate
synthetic surface velocity data using the
``true'' logarithmic basal sliding field,
$\beta_{\text{true}}\left(x\right):=
\log\left(1\,200+1\,100\sin\left(\frac{2\pi
	x}{W}\right)\right)$. Given this basal sliding field, we solve
Equations~\ref{Stokeseqn:1}-\ref{Stokeseqn:4}, 
extract the tangential velocity component at $100$
uniformly distributed points on the top boundary $\Gamma_t$, and add
$1\%$ relative Gaussian noise to each data point, resulting in the
synthetic data $\bm{d}$.

It remains to define the prior distribution for the parameter field
$\beta$.  The average value of $\beta_{\text{true}}$ is used as
constant prior mean $\overline{\beta}\left(x\right)= 6.73315 \approx
\frac{1}{W}\int_{0}^{W}\beta_{\text{true}}\left(s\right)\mathrm{d}s$. The
prior covariance matrix $\bm{\Gamma}_{\text{prior}}$ is a
discretization of the covariance PDE operator $\mathcal{C}:=\left(\delta
I-\gamma\Delta\right)^{-1}$, with $\gamma=6\times 10^{2}$ and $\delta
=2.4\times 10^{-3}$, with Robin boundary conditions~\cite{daonstadler2018}. These values are chosen in order to provide a
relatively large prior correlation length of
$10^{3}\left[\text{m}\right]$~\cite{lindgrenruelindstrom2011}.
Next, we summarize the computation of the MAP
point and the compression of the Gauss-Newton data-misfit Hessian matrix at the MAP point.

\pgfmathsetseed{2}

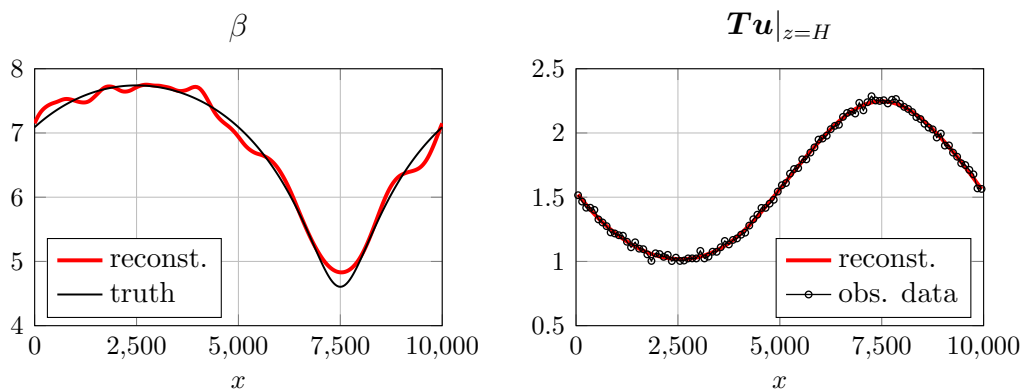
\begin{figure}[tb] 
	\begin{center}
		\begin{tabular}{rl}
			\begin{tikzpicture}[baseline, trim axis left]
			\begin{axis}[
			grid=major,
			width=7cm,
			height=5cm,
			xlabel = $x$,
			xmin=0.0, xmax = 1.e4,
			ymin = 4.0, ymax = 8.0,
			ytick={4.0, 5.0, 6.0, 7.0, 8.0},
			xtick={0.0, 2.5e3, 5.0e3, 7.5e3, 1.0e4},
			legend style={font=\small, nodes=right},
			legend pos=south west,
			title=$\beta$
			]
			\addlegendentry{reconst.}
			\addplot [color=red, line width = 1.5pt]
			table {2DbetaReconstruction.dat};
			\addlegendentry{truth}
			\addplot [color=black, line width = 0.75pt]
			table {./2DbetaTruth.dat};
			\end{axis}
			\end{tikzpicture}
			&
			\begin{tikzpicture}[baseline, trim axis right]
			\begin{axis}[
			grid=major,
			width=7cm,
			height=5cm,
			xlabel = $x$,
			xmin=0.0, xmax = 1.e4,
			ymin = 0.5, ymax = 2.5,
			ytick={0.5, 1.0, 1.5, 2.0, 2.5},
			xtick={0.0, 2.5e3, 5.0e3, 7.5e3, 1.0e4},
			legend style={font=\small, nodes=right},
			legend pos=south east,
			title = $\bm{T}\bm{u}|_{z=H}$
			]
			\addlegendentry{reconst.}
			\addplot [color=red, line width = 1.5pt]
			table {2DuxReconstruction.dat};
			\addlegendentry{obs. data}
			\addplot [color=black, line width = 0.5pt, mark=o, mark size=1.25pt]
			table {2DuxObserved.dat};
			\end{axis}
			\end{tikzpicture}\\
		\end{tabular}
	\end{center}
	\caption{Shown for Example I are on the left the MAP point $\beta^{\star}$ (red) and the
		truth basal sliding parameter $\beta_{\text{true}}$ (black) used to
		generate synthetic observations of the tangential velocity
		component on the upper surface $\Gamma_{t}$. Shown on the right are noisy synthetic
		observations (black dots) used for computing the MAP point  and the
		associated tangential surface velocity reconstruction (red).}
	\label{fig:ISMIP:MAP}
\end{figure} 

\subsection{MAP point and HODLR Gaussianized posterior}
The nonlinear optimization problem for finding the MAP estimate is
solved using an inexact Gauss-Newton minimization method with backtracking linesearch~\cite{nocedalwright2006}, where
the linear systems are iteratively solved by the conjungate
gradient method. The resulting MAP point is
shown in Figure~\ref{fig:ISMIP:MAP}. The MAP parameter field $\beta^{\star}$, closely resembles the true parameter
$\beta_{\text{true}}$, which is a consequence of the large amount of
available data and relatively small noise level.

\begin{figure}[tb]
	\begin{center}
		\includegraphics[scale=.35]{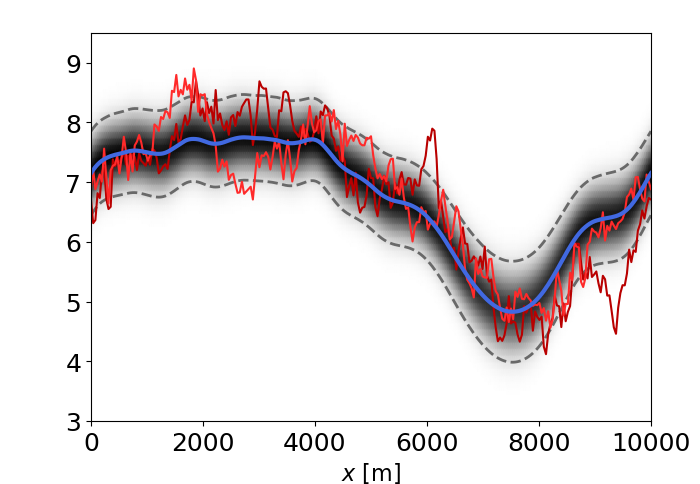} 
		\hfil 
		\includegraphics[scale=.35]{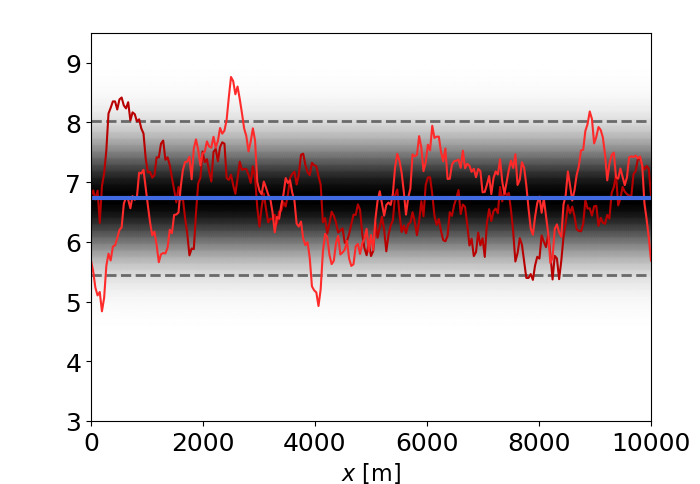}\hfil
	\end{center}
	\caption{Results for Example I: Two random samples (red), mean $\overline{\beta}$ (blue)
		and boundaries of the region $R=\lbrace \left(x,y\right) \text{ such
			that } 0\leq x\leq W \text{ and
		}\overline{\beta}(x)-2\sigma(x)\leq y\leq
		\overline{\beta}(x)+2\sigma(x)\rbrace$ (dashed black) are shown for
		the prior (left) and a HODLR Gaussianized posterior using the scheme
		described in Section~\ref{subsec:HODLRGaussianizedPosterior}
		(right).}\label{fig:ISMIP:samples}
\end{figure}

Next, we use the Gaussianized posterior distribution with a compressed
prior-preconditioned data-misfit 
Hessian $\bm{H}_{\text{misfit}}^{\prime}$ to generate approximate samples from the posterior
distribution.  Upon construction of the HODLR compression of the prior-preconditioned data-misfit
Hessian (details and comparisons can be found below in Section~\ref{subsec:ISMIP:props}), we draw samples from the HODLR Gaussianized
posterior as outlined in
Section~\ref{subsec:HODLRGaussianizedPosterior}. 
In Figure~\ref{fig:ISMIP:samples}, we compare the mean, pointwise standard
deviation and samples from the prior and the posterior distributions.
As expected, we find that the data updates our
belief about the spatially distributed parameter field and reduces the
uncertainty. In particular, the $2\sigma$ bounds on the
one-dimensional point marginals $\sigma\left(x\right)$,
$\bm{\sigma}_{i}=\left[\bm{\Gamma}_{i,i}\right]^{-1/2}$ of the
Gaussianized posterior and the prior distributions are shown, in order
to verify that the samples are largely contained within two standard
deviations of their respective means. The prior-preconditioned data-misfit
Hessian $\bm{H}_{\text{misfit}}^{\prime}$, is compressed using a relative tolerance of $10^{-6}$, that is~$\|\bm{H}_{\text{misfit}}^{\prime}-\bm{\tilde{H}}_{\text{misfit}}^{\prime}\|_{2}/\|\bm{H}_{\text{misfit}}^{\prime}\|_{2}\leq 10^{-6}$,
with high probability.

\subsection{Dependence of Hessian block spectra on problem setting}\label{subsec:ISMIP:props}
Next, we study how problem features impact the
numerical suitability of using global low-rank and HODLR
compressions to approximate the Gauss-Newton data-misfit Hessian.
 In this and
subsequent sections we measure the cost to generate the matrix compression
in terms of Hessian vector products, which we also describe as Hessian applies, as each said vector product
requires two linearized PDE solves and thus dominates the
computational cost. We use the result of Appendix~\ref{subsec:erroranalysis}, 
to claim $\varepsilon$ absolute error in a level $L$ HODLR approximation, when there is no more than $\varepsilon/L$ absolute error in each off-diagonal block.
 What is particular to this section, is that \textit{adaptive} single-pass and HODLR algorithms are used to generate global low-rank and HODLR approximations, based on absolute tolerance criteria. The absolute tolerance algorithmic input is scaled by the largest global low-rank singular value in order to report relative approximation errors. We note that additional errors are neglected in the reported approximation error such as that incurred in the peeling process~\cite{linlulexing2011, martinsson2016} and additional approximation assumptions in the single-pass algorithm, both of which are not expected to be significant.

\begin{figure}[tb]
	\begin{center} 
		\begin{tikzpicture}[baseline, trim axis right]
		\begin{axis}[
		grid=major,
		width=8cm,
		height=6cm,
		xlabel = $\|\bm{H}_{\text{misfit}}^{\text{GN}}-\bm{\tilde{H}}_{\text{misfit}}^{\text{GN}}\|_{2}/\|\bm{H}_{\text{misfit}}^{\text{GN}}\|_{2}\text{, approximation error}$,
		ylabel = {Computaitonal cost (Hessian applies)},
		xmin=1.e-8, xmax = 1.e-2,
		xmode=log,
		ymin = 0.0, ymax = 125.0,
		xtick={1.e-8, 1.e-6, 1.e-4, 1.e-2, 1.0},
		ytick={0.0, 25.0, 50.0, 75.0, 100.0, 125.0},
		legend style={font=\small, nodes=left,
			nodes={scale=0.7, transform shape}},
		legend pos=outer north east
		]
		\addlegendentry{HODLR, $\phi=1/200$}
		\addplot [color=blue, line width = 1.25pt]
		table {2DAR200HODLRcompression.dat};
		\addlegendentry{HODLR, $\phi=1/100$}
		\addplot [color=blue, dashdotted, line width = 1.25pt]
		table {2DAR100HODLRcompression.dat};
		\addlegendentry{HODLR, $\phi=1/50$\phantom{e}}
		\addplot [color=blue, dashed,line width = 1.25pt]
		table {2DAR50HODLRcompression.dat};
		\addlegendentry{HODLR, $\phi=1/25$\phantom{e}}
		\addplot [color=blue, loosely dotted,line width = 1.25pt]
		table {2DAR25HODLRcompression.dat};
		\addlegendentry{LR, $\phi=1/200$}
		\addplot [color=black, line width = 1.25pt]
		table {2DAR200LRcompression.dat};
		\addlegendentry{LR, $\phi=1/100$}
		\addplot [color=black, dashdotted, line width = 1.25pt]
		table {2DAR100LRcompression.dat};
		\addlegendentry{LR, $\phi=1/50$\phantom{e}}
		\addplot [color=black, dashed, line width = 1.25pt]
		table {2DAR50LRcompression.dat};
		\addlegendentry{LR, $\phi=1/25$\phantom{e}}
		\addplot [color=black, loosely dotted, line width = 1.25pt]
		table {2DAR25LRcompression.dat};
		\end{axis}
		\end{tikzpicture}
		\caption{Comparison of HODLR and global low-rank (LR) compression costs of the
			Gauss-Newton data-misfit Hessian $\bm{H}_{\text{misfit}}^{\text{GN}}$,
			for Example I with ice sheet aspect ratio $\phi$. This figure shows that for
			low aspect ratios, HODLR becomes more efficient than global low-rank
			for medium levels of target accuracy.}
		\label{fig:aspectratiostudy}
	\end{center}
\end{figure}
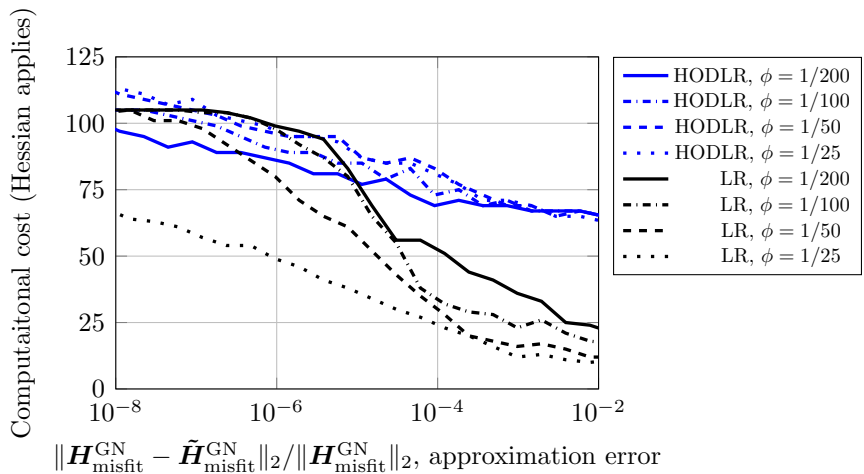

\begin{figure}[tb]
	\begin{center}
		\begin{tikzpicture}[baseline, trim axis right]
		\begin{axis}[
		grid=major,
		width=8cm,
		height=6cm,
		xlabel = $\|\bm{H}_{\text{misfit}}^{\text{GN}}-\bm{\tilde{H}}_{\text{misfit}}^{\text{GN}}\|_{2}/\|\bm{H}_{\text{misfit}}^{\text{GN}}\|_{2}\text{, approximation error}$,
		ylabel = Computaitonal cost (Hessian applies),
		xmin=1.e-8, xmax = 1.e-2,
		xmode=log,
		ymin = 0.0, ymax = 150.0,
		xtick={1.e-8, 1.e-6, 1.e-4, 1.e-2, 1.0},
		ytick={0.0, 25.0, 50.0, 75.0, 100.0, 125.0, 150.0, 
			175.0},
		legend style={font=\small, nodes=left, nodes={scale=0.7, transform shape}},
		legend pos=outer north east,
		]
		\addlegendentry{HODLR, $\text{dim}\left(\bm{\beta}\right)=128$}
		\addplot [color=blue, dashdotted, line width = 1.25pt]
		table {2Dn128HODLRcompression.dat};
		\addlegendentry{HODLR, $\text{dim}\left(\bm{\beta}\right)=256$}
		\addplot [color=blue, dashed, line width = 1.25pt]
		table {2Dn256HODLRcompression.dat};
		\addlegendentry{HODLR, $\text{dim}\left(\bm{\beta}\right)=512$}
		\addplot [color=blue, densely dotted, line width = 1.25pt]
		table {2Dn512HODLRcompression.dat};
		\addlegendentry{LR, $\text{dim}\left(\bm{\beta}\right)=128$}
		\addplot [color=black, dashdotted, line width = 1.25pt]
		table {2Dn128LRcompression.dat};
		\addlegendentry{LR, $\text{dim}\left(\bm{\beta}\right)=256$}
		\addplot [color=black, dashed, line width = 1.25pt]
		table {2Dn256LRcompression.dat};
		\addlegendentry{LR, $\text{dim}\left(\bm{\beta}\right)=512$}
		\addplot [color=black, densely dotted, line width = 1.25pt]
		table {2Dn512LRcompression.dat};
		\end{axis}
		\end{tikzpicture}
	\end{center}
	\caption{Dependence of HODLR and global low-rank (LR) compression costs of the
		Gauss-Newton data-misfit Hessian on $\text{dim}\left(\bm{\beta}\right)$, the
		dimension of the discretized logarithmic basal sliding field
		for Example I. The cost of global low-rank compression is almost constant, while the cost of HODLR compression increases as the mesh is refined.}
	\label{fig:paramdimstudy}
\end{figure}
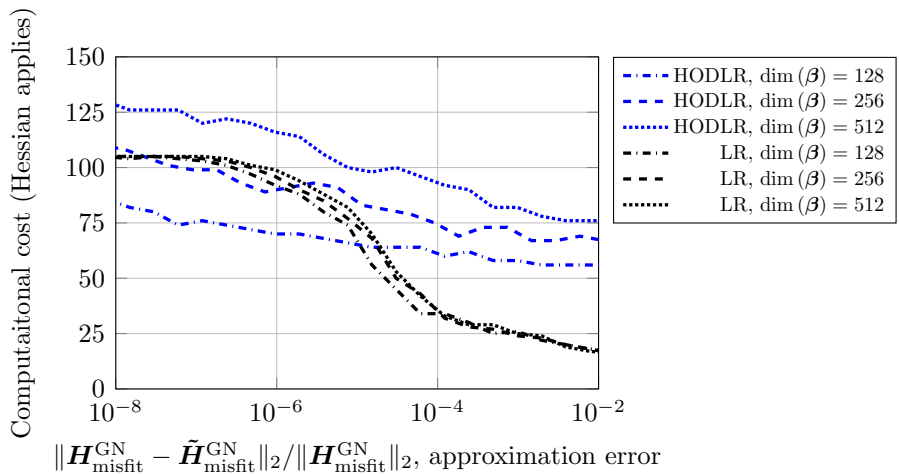

\paragraph{Influence of aspect ratio}
Here, we vary the aspect ratio of the domain $\phi=H/W$, where $H$ and $W$ are the domain height and width respectively,
in order to study how it
influences the block spectra of the Gauss-Newton data-misfit Hessian and
ultimately the computational cost. Figure~\ref{fig:aspectratiostudy}
shows that the global spectrum is more sensitive to changes in the relative
length scale $\phi$ than the spectra of the off-diagonal blocks.  Low-rank
approximations of the off-diagonal blocks become computationally
cheaper as $\phi$ decreases as a result of the sensitivity cones
becoming increasingly localized as the ice sheet thickness decreases. Global low-rank approximations become more expensive
as $\phi$ decreases, a result of the data being more informative. We
note that realistic problems, such as the Humboldt glacier and the
Greenland ice sheet studied later in Section \ref{sec:HumboldtGreenland},
have small aspect ratios and are thus expected to have data-misfit Hessians that
are less amenable to global low-rank approximation.
\paragraph{Influence of the parameter dimension}
We now vary the level of mesh discretization refinement in order to study the influence of data informativeness, through the discretized parameter dimension $N=\text{dim}(\bm{\beta})$, on the computational cost to generate
HODLR and global low-rank approximations of the Gauss-Newton data-misfit Hessian. The hierarchical depth $L$ is
incremented for every doubling of the discretized parameter dimension,
in order that the hierarchical depth scales with the logarithm of the size of the Hessian matrix, 
a condition described in Section~\ref{subsec:HODLRdef}. 
Figure~\ref{fig:paramdimstudy} provides computational evidence of the
claim made in Section~\ref{subsec:HODLRdef}, that the number of
applies needed to hierarchically compress an operator with HODLR
structure is $\mathcal{O}\left(\log\,N\right)$. On the contrary, the number
of applies to generate the global low-rank approximation is rather insensitive
to the level of mesh refinement. 

\begin{figure}[tb]
	\begin{center}
		\begin{tikzpicture}[baseline, trim axis right]
		\begin{axis}[
		grid=major,
		width=8cm,
		height=6cm,
		xlabel = $\|\bm{H}_{\text{misfit}}^{\text{GN}}-\bm{\tilde{H}}_{\text{misfit}}^{\text{GN}}\|_{2}/\|\bm{H}_{\text{misfit}}^{\text{GN}}\|_{2}\text{, approximation error}$,
		ylabel = Computational cost (Hessian applies),
		xmin=1.e-8, xmax = 1.e-2,
		ymin = 0.0, ymax = 200.0,
		xmode=log,
		ytick={0.0, 25.0, 50.0, 75.0, 100.0, 125.0, 150., 
			175.0, 200.0},
		legend style={font=\small, nodes=left, nodes={scale=0.7, transform shape}},
		legend pos=outer north east,
		]
		\addlegendentry{HODLR, $\text{dim}(\bm{d})=1.0\times 10^{2}$}
		\addplot [color=blue, dashed, line width = 1.25pt]
		table {2Dnobs100HODLRcompression.dat};
		\addlegendentry{HODLR, $\text{dim}(\bm{d})=1.5\times 10^{2}$}
        \addplot [color=blue, dashdotted, line width = 1.25pt]
        table {2Dnobs150HODLRcompression.dat};
		\addlegendentry{HODLR, $\text{dim}(\bm{d})=2.0\times 10^{2}$}
		\addplot [color=blue, densely dotted, line width = 1.25pt]
		table {2Dnobs200HODLRcompression.dat};
		\addlegendentry{LR, $\text{dim}(\bm{d})=1.0\times 10^{2}$}
		\addplot [color=black, dashed, line width = 1.25pt]
		table {2Dnobs100LRcompression.dat};
		\addlegendentry{LR, $\text{dim}(\bm{d})=1.5\times 10^{2}$}
        \addplot [color=black, dashdotted, line width = 1.25pt]
        table {2Dnobs150LRcompression.dat};
		\addlegendentry{LR, $\text{dim}(\bm{d})=2.0\times 10^{2}$}
		\addplot [color=black, densely dotted, line width = 1.25pt]
		table {2Dnobs200LRcompression.dat};
		\end{axis}
		\end{tikzpicture}
		\caption{Dependence of HODLR and global low-rank (LR) compression costs of the
			Gauss-Newton data-misfit Hessian on $\text{dim}(\bm{d})$, the data dimension, for Example I. The computational cost
			for global low-rank approximation increases with
			the number of observations, while the cost for HODLR compression is
			rather insensitive.}
		\label{fig:nobsstudy}
	\end{center}
\end{figure}
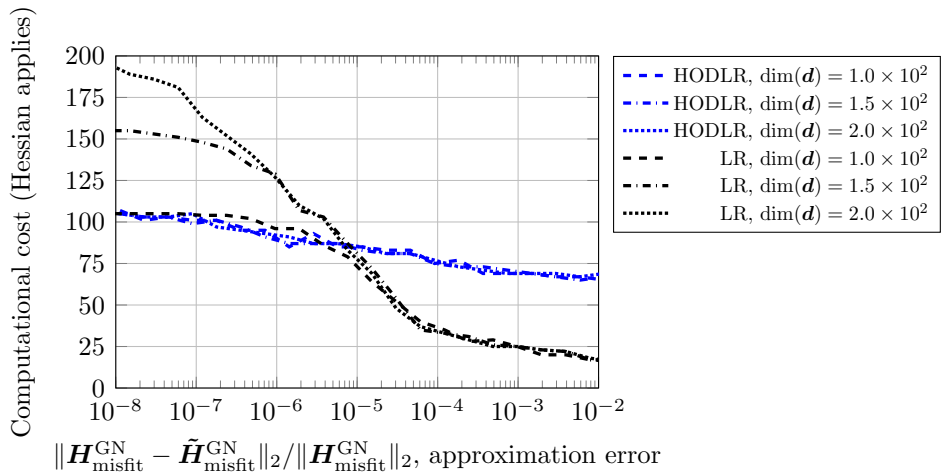

\paragraph{Influence of the data dimension}
Figure~\ref{fig:nobsstudy} shows that the global rank grows with the
number of observations points and thus global low-rank compression tends to be
less efficient for problems with strongly informative observation data.  The rate of
spectral decay of the (Gauss-Newton) data-misfit Hessian is related to the degree of ill-posedness
of the unregularized inverse problem.
As the number of observations increases, these associated model predictions  are increasingly sensitive to small scale variations in the basal sliding
field. Thus, more data, generally makes the data set more informative about the
parameter and the (Gauss-Newton) data-misfit Hessian have a weaker rate of spectral decay.

\section{Example II: Humboldt glacier and Greenland ice sheet}
\label{sec:HumboldtGreenland}
Here, we study the scalability of the proposed methods using
large-scale ice sheet problems which are typically used in climate
simulations. Namely, we focus on the Humboldt glacier in North-West
Greenland, and the entire Greenland ice sheet. For these simulations,
we use the ice sheet model MALI,~\cite{hoffman2018}, which relies on
Albany,~\cite{tezaurperegoetal2015}, a C++ multi-physics library for
the implementation of the so-called first-order approximation of
Stokes equations. This first-order approximation is based on scaling arguments
motivated by the shallow nature of ice sheets and uses the
incompressibility condition to reduce the unknows to the horizontal
velocities. We use PyAlbany~\cite{liegeois2022} a convenient Python
interface to the Albany package, which in turn builds upon 
Trilinos~\cite{trilinos-website}. Albany is designed to support parallel and
scalable finite-element discretized PDE solvers and various analysis
capabilities. Details about the parameter, state, data dimensions as
well as the number of cores and hierarchical levels used in the
computations is provided in Table~\ref{table:exampleII}.
\begin{table}[htp] 
	\centering 
	\begin{tabular}{|c|c|c|}
		\hline
		& Humboldt & Greenland \\
		\hline 
		dim$(\bm{\beta})$ & $11\,608$ & $320\,116$ \\
		\hline
		dim$(\bm{u})$ & $255\,376$ & $7\,042\,552$ \\
		\hline
		dim$(\bm{d})$ & $23\,216$ & $640\,232$ \\
		\hline
		\# of cores & $120$ & $2\,048$ \\
		\hline 
		$L$ & $8$ & $10$ \\
		\hline
	\end{tabular}
	\caption{Problem specifications for the Humboldt glacier and Greenland
		ice-sheet problems (Example II). Dimension of the discretized parameter field
		dim$(\bm{\beta}$), dimension of the discretized velocity field
		dim$(\bm{u})$, dimension of the observations
		dim$(\bm{d})$, processors employed for computations and
		$L$ depth of HODLR hierarchical partitioning.}
	\label{table:exampleII}
\end{table}

The following study is partially motivated by findings made in the
Section~\ref{sec:2DResults}, namely that the role of the aspect ratio between the
vertical and horizontal directions (see Section~\ref{subsec:ISMIP:props}) influences the ability to use global low-rank
compression and favors HODLR compression. We generate HODLR and global low-rank approximations
and then based on the computed spectra, Equation~\ref{eq:HODLRcost} and $\zeta^{\text{LR}}=r+d$, we estimate
the computational cost. Additionally, we study
how the ordering of the degrees of freedom impacts the spectral decay
for off-diagonal blocks of the data-misfit Hessian.
We present
results for both, the Humboldt glacier, which expands about $4\times
10^{2}$ [km] laterally, and the Greenland ice sheet, which
expands about $1.8\times 10^{3}$ [km]. The ice is at most~$3.4$~[km] thick, resulting in approximate
aspect ratios of $8.5\times 10^{-3}$ for Humboldt and $1.9\times 10^{-3}$ for
Greenland. We use a nonuniform triangulation of the Greenland ice
sheet, with mesh size ranging from 1 to 10 [km], and we then extrude
it in the vertical direction, obtaining a 3D mesh having 10 layers of
prismatic elements. The velocity observations at the top surface of the Greenland ice sheet are obtained from satellite observations~\cite{joughin2015}. The MAP basal sliding field and the temperature fields are obtained as part of the initialization process, using a numerical optimization approach to match the ice velocity observations and constrained by the first-order flow model coupled
with a temperature model~\cite{perego2022}.
Additional details about the mesh geometries and data, in particular regarding the Humboldt glacier, can be found in~\cite{hillebrand2022}.

\begin{figure}[tb]
	\begin{center}
		\includegraphics[scale=.1125]{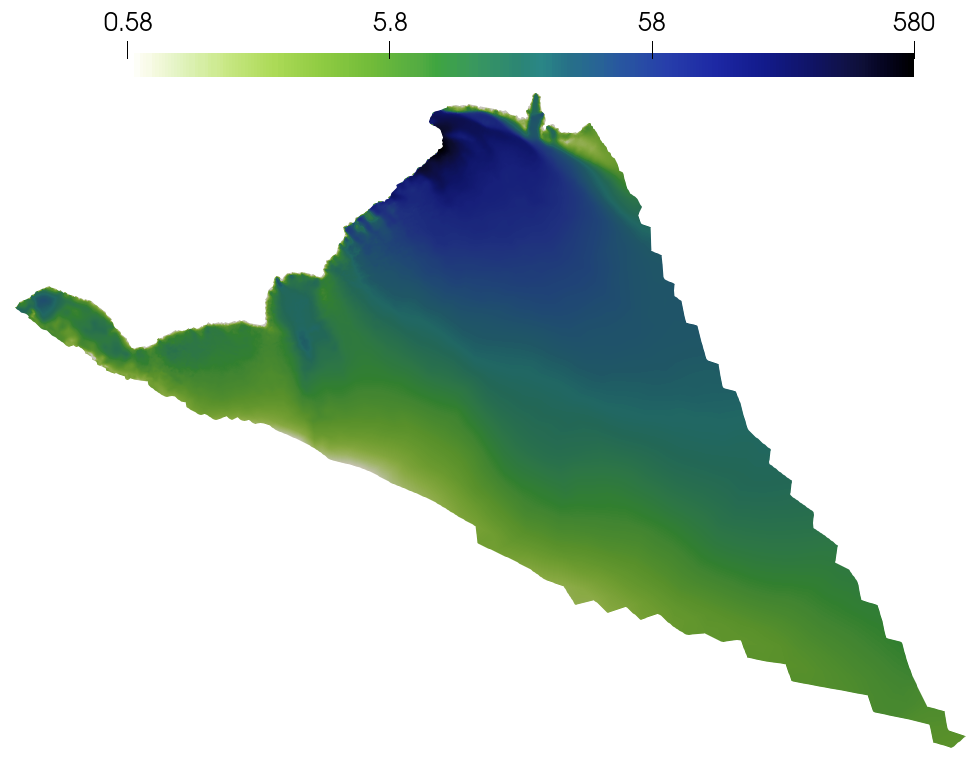} \hfil 
		\includegraphics[scale=.1125]{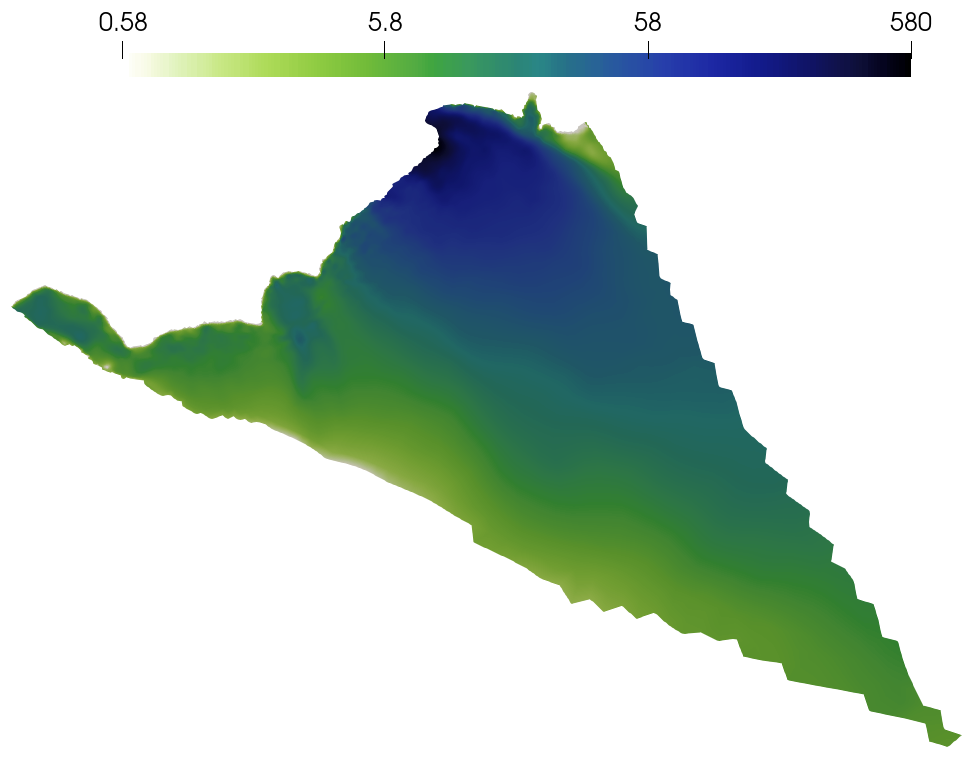} \hfil 
		\includegraphics[scale=.1125]{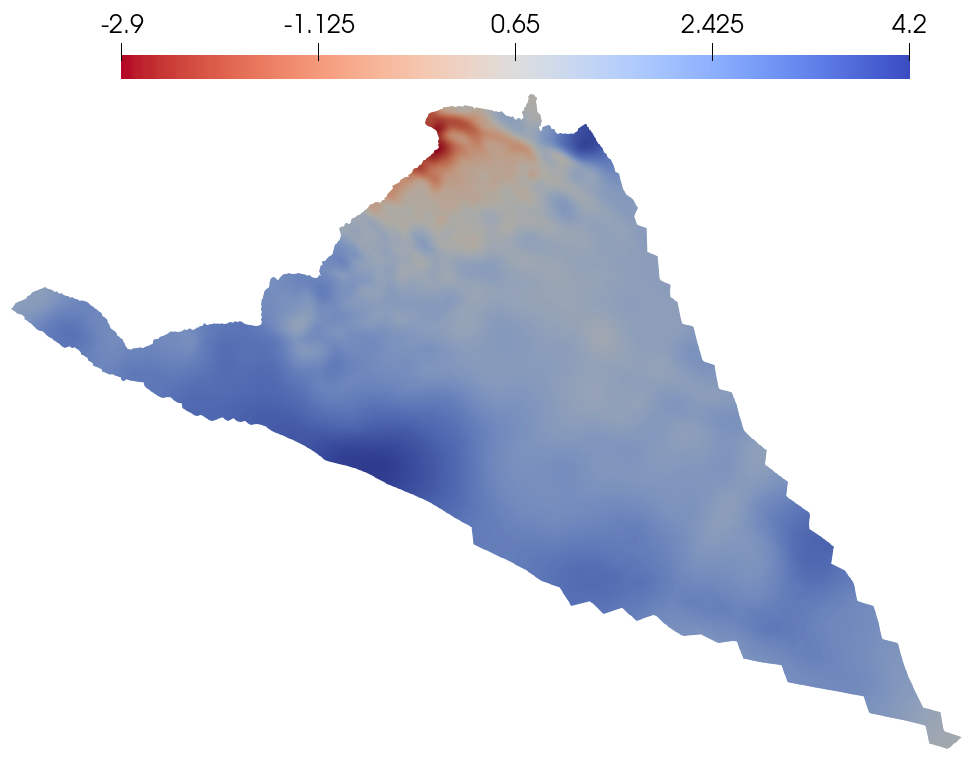} \hfil \\
		\includegraphics[scale=.1625]{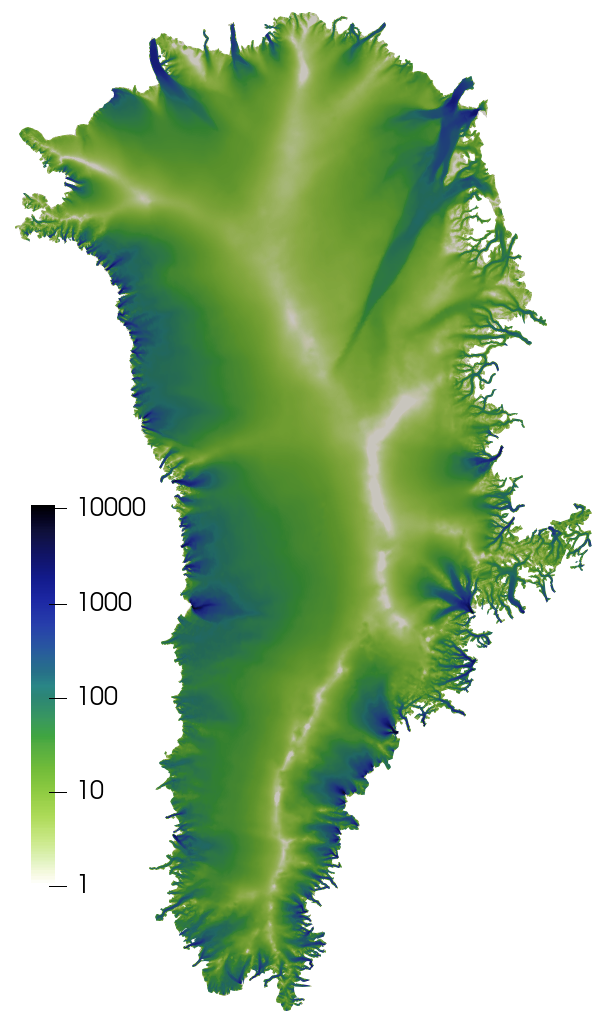} \hfil 
		\includegraphics[scale=.1625]{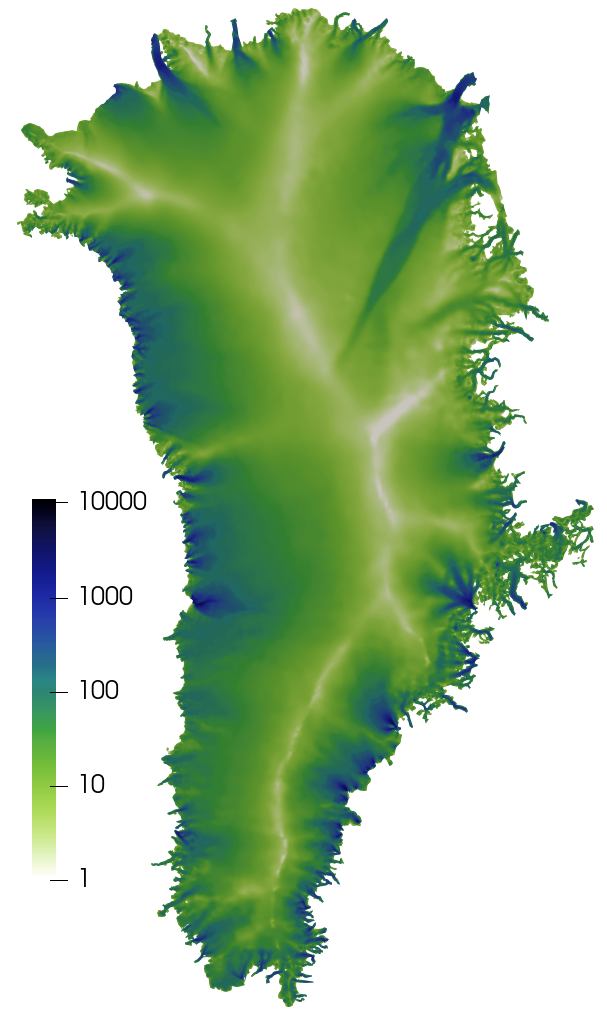}	\hfil 
		\includegraphics[scale=.1625]{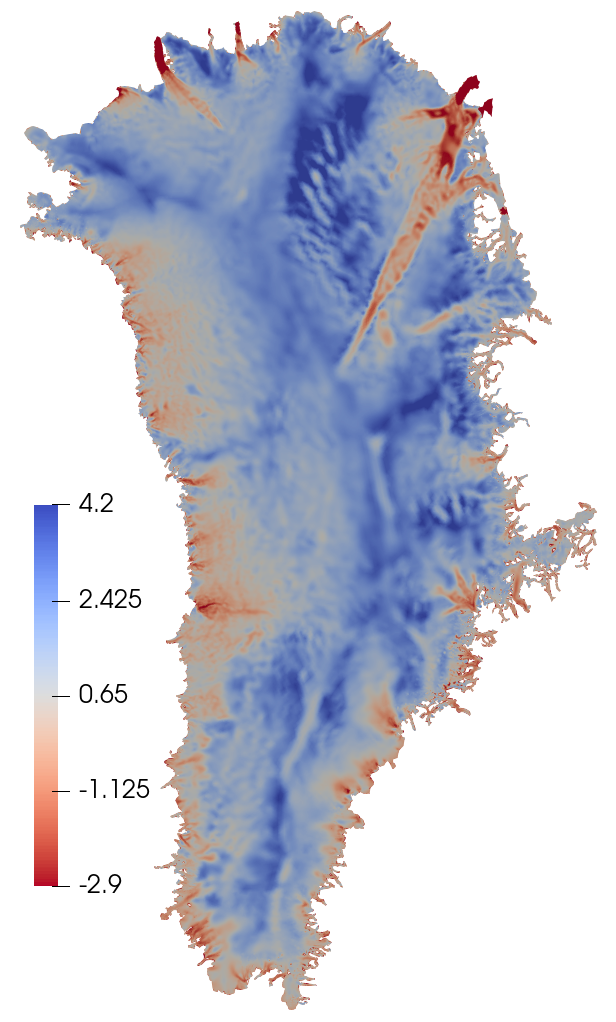} \hfil
	\end{center} 
	\caption{Data and MAP estimates for Example II. Shown are the surface velocity observation
		data (left), and the reconstructed
		surface velocity field (middle) that is based on the MAP
		estimate of the logarithmic basal sliding field
		(right). Top row is for the Humboldt glacier and bottom row for the
		Greenland ice sheet.}
	\label{fig:Humboldtfieldplots}
\end{figure}

In Figure~\ref{fig:Humboldtfieldplots}, we show the observed surface
velocity $\bm{d}$ in [m/yr], the MAP estimates of the logarithmic
basal sliding field $\beta^{\star}$ ($\exp(\beta^{\star})$ is in
[kPa yr/m]) and surface velocity in [m/yr] generated by the model.

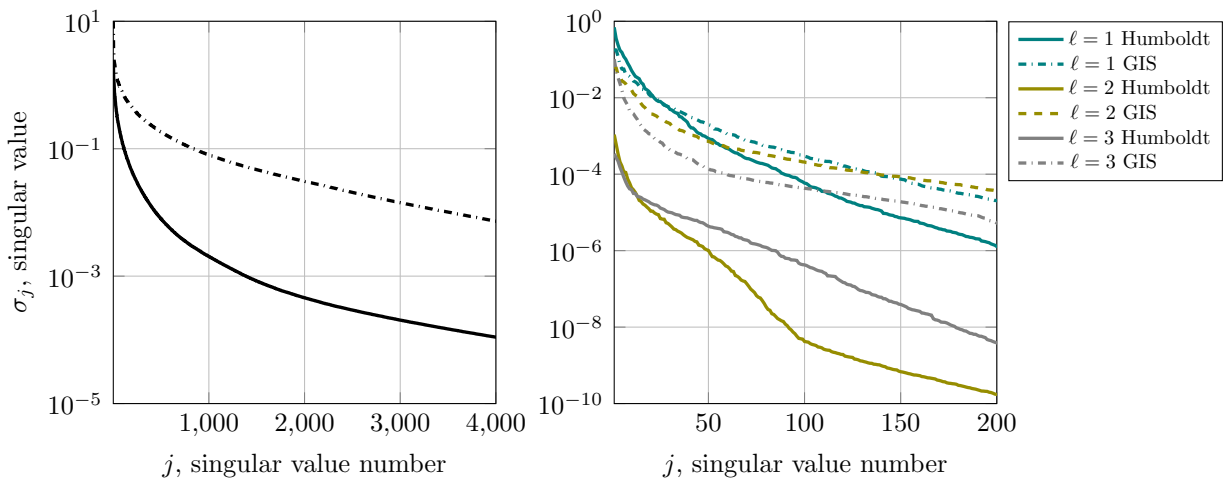
\begin{figure}[tb] 
	\begin{center}
		\begin{tabular}{rl}
			\begin{tikzpicture}[baseline, trim axis right]
			\begin{axis}[
			grid=major, width=0.425\textwidth, height=0.425\textwidth,
			ymode=log,
			xlabel = {$j$, singular value number},
			ylabel = {$\sigma_{j}$, singular value},
			xmin=1, xmax = 4.e3,
			ymin=1.e-5, ymax = 1.e1,
			xtick={1.e3, 2.e3, 3.e3, 4.e3},
			ytick={1.e-5, 1.e-3, 1.e-1, 1.e1},
			]
			\addplot [color=black, line width = 1.25pt]
			table{HumboldtSigGlb.dat};
			\addplot [color=black, dashdotted, line width = 1.25pt]
			table{GISSigGlb.dat};
			\end{axis}	
			\end{tikzpicture}
			&		
			\begin{tikzpicture}[baseline, trim axis right]
			\begin{axis}[
			grid=major, width=0.425\textwidth, height=0.425\textwidth,
			ymode=log,
			xlabel = {$j$, singular value number},
			scaled ticks=false,
			xmin=1,xmax=200,
			ymin=1.e-10,ymax=1.e0,
			ytick={1.e-10, 1.e-8, 1.e-6, 1.e-4, 1.e-2, 1.e0},
			legend style={font=\small, nodes=left, nodes={scale=0.7, transform shape}},
			legend pos=outer north east,
			]
			\addlegendentry{$\ell=1$ Humboldt}
			\addplot [color=teal, line width = 1.25pt]
			table{HumboldtSigL0J0.dat};
			\addlegendentry{$\ell=1$ GIS\phantom{boldt }}
			\addplot [color=teal, dashdotted, line width = 1.25pt]
			table{GISSigL0J0.dat};
			\addlegendentry{$\ell=2$ Humboldt}
			\addplot [color=olive, line width = 1.25pt]
			table{HumboldtSigL1J0.dat};
			\addlegendentry{$\ell=2$ GIS\phantom{boldt }}
			\addplot [color=olive, dashed, line width = 1.25pt]
			table{GISSigL1J0.dat};
			\addlegendentry{$\ell=3$ Humboldt}
			\addplot [color=gray, line width = 1.25pt]
			table{HumboldtSigL2J0.dat};
			\addlegendentry{$\ell=3$ GIS\phantom{boldt }}
			\addplot [color=gray, dashdotted, line width = 1.25pt]
			table{GISSigL2J0.dat};
			\end{axis}
			\end{tikzpicture}
		\end{tabular}		
	\end{center}
	\caption{Singular values of the data-misfit Hessian (left
		figure) and various off-diagonal blocks of the data-misfit
		Hessian (right figure) for Example II. The color-scheme in
		the right most figure is consistent with Figure~\ref{fig:hmatrixpartitioningstructure}. On the left, the
		singular values of the Humboldt and Greenland data-misfit
		Hessians are shown using a solid and dash-dotted line,
		respectively. On the right, we show the singular values of
		the upper most
		blocks, that is $\bm{A}^{\left(\ell\right)}_{1,2}$ as defined
		in Section~\ref{subsec:erroranalysis}.}
	\label{fig:HumboldtGreenlandSpectra}
\end{figure}

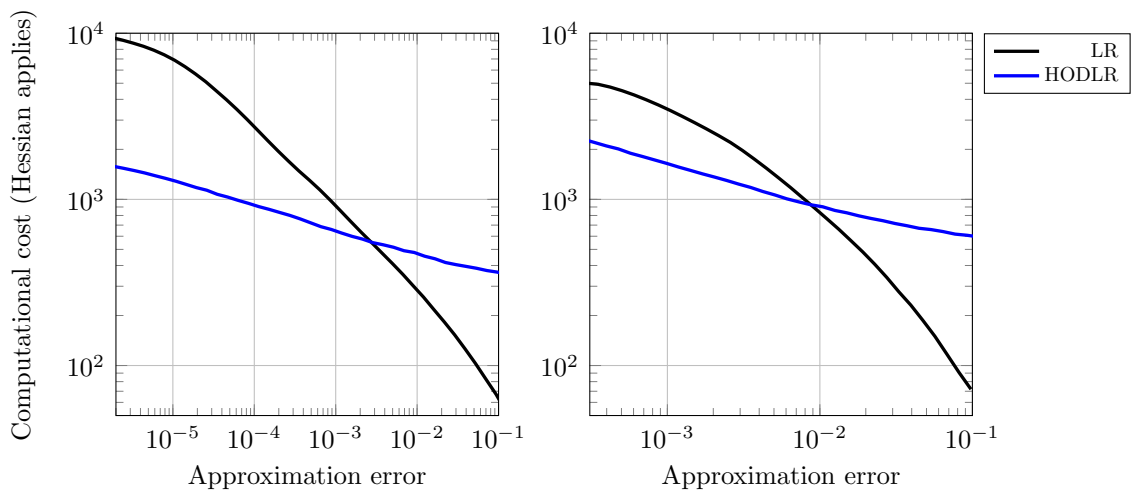
\begin{figure}[tb] 
	\begin{center}
		\begin{tabular}{rl}
			\begin{tikzpicture}[baseline, trim axis right]
			\begin{axis}[
			grid=major, width=0.425\textwidth, height=0.425\textwidth,
			xmode=log, ymode=log,
			xlabel = {Approximation error},
			ylabel = {Computational cost (Hessian applies)},
			xmin=2.e-6, xmax = 1.e-1,
			ymin=5.e1, ymax = 1.e4,
			legend style={font=\small, nodes=left, nodes={scale=0.7, transform shape}},
			legend pos=outer north east,
			]
			\addplot [color=black, line width = 1.25pt]
			table{HumboldtLRcompression.dat};
			\addplot [color=blue, line width = 1.25pt]
			table{HumboldtHODLRcompression.dat};
			\end{axis}
			\end{tikzpicture}
			&
			\begin{tikzpicture}[baseline, trim axis right]
			\begin{axis}[
			grid=major, width=0.425\textwidth, height=0.425\textwidth,
			xmode=log, ymode=log,
			xlabel = {Approximation error},
			xmin=3.1e-4, xmax = 1.e-1,
			ymin=5.e1, ymax = 1.e4,
			legend style={font=\small, nodes=left, nodes={scale=0.7, transform shape}},
			legend pos=outer north east,
			]
			\addlegendentry{LR}
			\addplot [color=black, line width = 1.25pt]
			table{GISLRcompression.dat};
			\addlegendentry{HODLR}
			\addplot [color=blue, line width = 1.25pt]
			table{GISHODLRcompression.dat};
			\end{axis}
			\end{tikzpicture}	
		\end{tabular} 
		\caption{Estimated computational costs (measured by the number of Hessian
			applies) to compress the Humboldt glacier (left) and
			Greenland ice-sheet (right) data-misfit Hessians into the
			global low-rank (LR) and hierarchical off-diagonal low-rank
			(HODLR) formats as a function of the approximation
			error $\|\bm{H}_{\text{misfit}}-\bm{\tilde{H}}_{\text{misfit}}\|_{2}/\|\bm{H}_{\text{misfit}}\|_{2}$.}
		\label{fig:GreenlandHumboldtError}
	\end{center}
\end{figure}

\subsection{HODLR compressability}
We next generate global low-rank approximations of a Greenland and Humboldt
data-misfit Hessian as well as low-rank approximations of various
off-diagonal blocks. Plots of the estimated singular values are
provided in Figure~\ref{fig:HumboldtGreenlandSpectra}. We observe that
the spectrum of the Greenland ice sheet decays substantially slower
than the one for the Humboldt glacier. Besides the different
sizes of these two discretized problems, this is also due to the
different aspect ratios. Having
estimated singular values of the data-misfit Hessians and the
appropriate off-diagonal blocks, one is able to estimate computational
costs to compress them into the global low-rank and HODLR matrix
formats. The computational cost as a function of Hessian approximation
target accuracy is given in Figure~\ref{fig:GreenlandHumboldtError},
wherein it is demonstrated that the HODLR compression format can offer
a favorable means to approximate data-misfit Hessians for large-scale
inverse problems governed by complex ice-sheet models.

\subsection{Impact of parameter degree of freedom ordering}
\label{subsec:dofordering}

We seek to ensure that the off-diagonal blocks,
determined by the hierarchical partitioning 
described in Section~\ref{subsec:HODLRdef},
of the data-misfit Hessian are low-rank.
For this reason, the nodes $\lbrace \bm{x}_{i}\rbrace_{i}$
associated to the degrees of freedom (dofs) are ordered
according to a kd-tree, i.e., a recursive hyperplane splitting.
The ordering provided by the kd-tree is such that
the $(i,j)$-entry of the distance matrix 
$\bm{D}_{i,j}=\|\bm{x}_{i}-\bm{x}_{j}\|_{2}$,
is typically small whenever $|i-j|$ is small,
that is the dof ordering preserves some notion of locality
(see Section~\ref{subsec:motivation}). 
In particular, a sparse permutation matrix
$\bm{B}$, is determined,
whose action reorders the dofs from the default
ordering provided by the finite element discretization to that specified
by the kd-tree. 
The data-misfit Hessian with respect to the kd-tree ordering, 
$\bm{H}^{\text{kd}}_{\text{misfit}}:=\bm{B}\bm{H}_{\text{misfit}}\bm{B}^{\top}$,
is then amenable to HODLR compression. Subsequently,
$\bm{B}^{\top}\bm{\tilde{H}}^{\text{kd}}_{\text{misfit}}\bm{B}$ is an approximation of the data-misfit Hessian with respect to the default ordering.

\begin{figure}[tb] 
	\begin{tikzpicture}
	\node at (0,7) {};
	\node at (0,-1) {};
	\node at (0,0) {
		\begin{axis}[compat=1.3,
		grid=major,
		width=10.5cm,
		height=7cm,
		ymode=log,
		xlabel = $j\text{, singular value number}$,
		ylabel = $\sigma_{j}\text{, singular value}$,
		xmin=1,xmax=500,
		ymin=1.e-11,ymax=1.e0,
		ytick={1.e-11, 1.e-9, 1.e-7, 1.e-5, 1.e-3, 1.e-1},
		legend style={nodes=left, nodes={scale=0.7, transform shape}},
		legend pos=north east,
		]
		\addlegendentry{default basis}
		\addplot [color=violet, line width = 1.pt]
		table{HumboldtkdtreeSigL0J0.dat};
		\addlegendentry{kd-tree basis}
		\addplot [color=violet, dash dot, line width = 1.5pt]
		table{HumboldtSigL0J0.dat};
		\end{axis}
	};
	\node at (10.1,2.0)
	{\includegraphics[scale=.2]{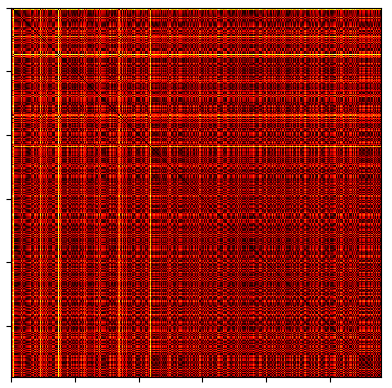}};
	\node at (2.3,1.4) {\includegraphics[scale=.2]{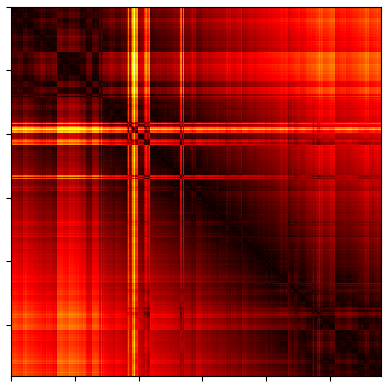}};
	\draw [-, opacity=0.6, very thick, white!50!black]
	(3.6,1.4)--(4.5,2.5);
	\draw [-, opacity=0.6, very thick, white!50!black] (8.8,2)--(7.5,3.9);
	\end{tikzpicture}
	\caption{Singular values of the hierarchical level $1$
		off-diagonal block, $\bm{A}_{1,2}^{(1)}$, of the Humboldt
		glacier data-misfit Hessian, when expressed in a kd-tree basis
		and the default basis.  Shown also are heat maps of the distance
		matrices $\bm{D}_{i,j}=\|\bm{x}_{i}-\bm{x}_{j}\|_{2}$, wherein
		the nodes $\lbrace \bm{x}_{i}\rbrace_{i}$, associated to the
		finite element degrees of freedom have been ordered according to
		a default standard and a kd-tree.
		\label{fig:humboldt-order}}
\end{figure}
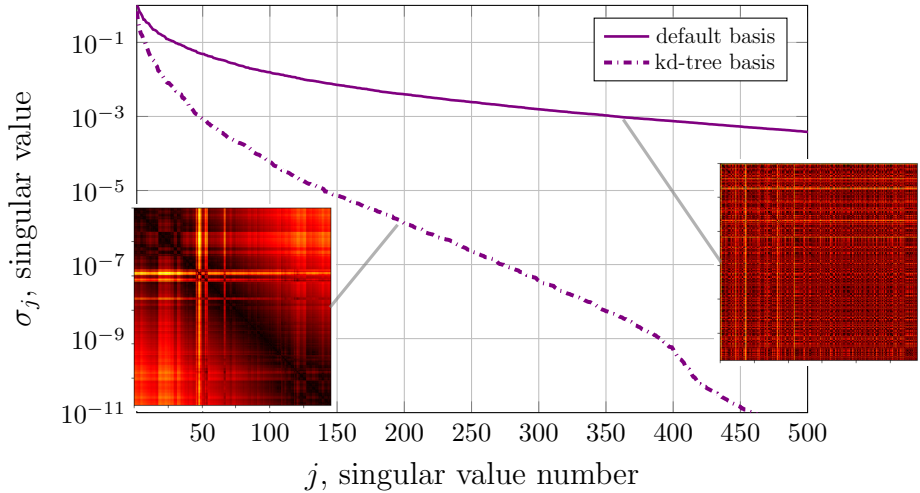

The dof ordering has no impact on a matrix's global numerical rank but
does indeed impact the numerical rank of its numerous submatrices that
are defined by a fixed partitioning scheme, such as the off-diagonal
blocks of an HODLR matrix (see Section~\ref{subsec:HODLRdef}).  Here,
we study the HODLR compressibility of the Humboldt glacier data-misfit
Hessian by comparing the rate of decay of an off-diagonal block's
singular values using the default ordering provided by Albany and the
ordering obtained by a kd-tree recursive hyperplane splitting.  As
observed in Figure~\ref{fig:humboldt-order}, the rate at which the
singular values of the level-1 off-diagonal block decay, strongly
depends on the dof ordering. This is because the ordering given by the
kd-tree better preserves locality, and as a consequence, by the
argument provided in Section~\ref{subsec:motivation}, the singular
values decay much faster when using the kd-tree ordering. The kd-tree
ordering therefore provides a substantially computationally cheaper
means to generate an HODLR approximation of the data-misfit
Hessian. Figure~\ref{fig:humboldt-order} also shows distance matrices
for the default and kd-tree bases. These show the improved locality
for the kd-orderings. Note that data-misfit Hessian matrices are
expected to follow a similar structure as these distance matrices,
which explains why the former's off-diagonal blocks can be compressed
more effectively in the kd-order than in the default order of dofs.
\section{Conclusion}

In this work, we motivated why data-misfit Hessians which arise from a class
of inverse problems governed by PDEs have HODLR matrix
structure. HODLR matrices can efficiently be inverted and factorized,
operations needed for solving inverse problems governed by PDEs by
Newton's method, for constructing Gaussian approximations and for
Markov chain Monte Carlo sampling methods. We study inverse ice sheet
problems, for which, under certain regimes, HODLR matrices provide a more computationally
efficient approximation format than the global low-rank matrix format. These
problems are those with highly informative data and small aspect ratio ice
sheets. While global low-rank matrices are favorable for large
discretized parameter dimension and small data dimension, we find that
HODLR matrices can offer computational savings for large-scale inverse problems
such as a Greenland ice sheet inverse problem with satellite
observational data and a discretized parameter dimension that exceeds
$10^{5}$.

For future work, we believe that the computational cost can be reduced further by utilizing hierarchical matrix partitionings
that satisfy a strong admissibility condition~\cite{hackbuschbohm2002},
as they are better suited to exploit data-misfit
Hessian structure. However, generating a hierarchical matrix
approximation with such a partitioning, e.g., by the
peeling method~\cite{linlulexing2011, martinsson2016},
requires substantially more Hessian vector products.
Ultimately, to further reduce the
computational cost of Hessian approximations in inverse problems
governed by PDEs, exploiting further problem
structure will be essential.

\section{Appendix}

\subsection{Randomized Compression Algorithms}
\label{subsec:randomizedcompressionalgorithms}

Here, for completeness we outline the randomized matrix-free double-pass global low-rank and HODLR compression algorithms. The essential ideas of the randomized double-pass low-rank algorithm~\cite{halkomartinsson2011} are
\begin{enumerate} 
	\item the application of a vector $\bm{\omega}$ with random entries to a matrix
	$\bm{A}$, yields a vector $\bm{y}=\bm{A}\bm{\omega}$, which is likely aligned with
	the dominant left singular vectors of $\bm{A}$;
	\item a matrix $\bm{Q}$, whose columns are nearly aligned with the dominant left singular vectors of $\bm{A}$, can be used to construct
	an accurate low-rank approximation $\bm{\tilde{A}}=\bm{Q}\bm{Q}^{\top}\bm{A}$ of $\bm{A}$.  
\end{enumerate} 
The double-pass randomized SVD algorithm is presented in Algorithm~\ref{alg:doublepass} and does not significantly differ from that
in~\cite{halkomartinsson2011}, specifically it is lines $7,8$ and $9$ that are distinct. This minor modification frees us from the need to compute a (parallel) singular
value decomposition (SVD) of a (distributed) $N\times k$ matrix, such
as $\bm{Z}$. Here, we only need to compute an SVD of the smaller
$k\times k$ matrix $\bm{R_{Z}}$. In the distributed memory parallelism
setting of Section~\ref{sec:HumboldtGreenland}, this algorithmic
modification allows us to only require the invocation of serial SVD routines, as 
$\bm{R_{Z}}$, which is typically small, is available on each processor.
\begin{algorithm}[H]
	\caption{Double-pass randomized SVD. \\
		\textbf{Input:} $\bm{A}\in\mathbb{R}^{N\times N}$, $r\in\mathbb{N}$ desired rank and oversampling parameter $d\in\mathbb{N}$.\\
		\textbf{Output:} low-rank approximation $\bm{\tilde{A}}$ of $\bm{A}$}
	\begin{algorithmic}[1]
		\STATE{$k=r+d$}
		\STATE{$\bm{\Omega}=$ \verb~randn~$(N,k)$ \hfill $\lbrace$Initiate random matrix$\rbrace$}
		\STATE{$\bm{Y}=\bm{A}\bm{\Omega}$ \hfill $\lbrace$Sample column space$\rbrace$}
		\STATE{$\bm{Q_{Y}}=$ \verb~orthog~($\bm{Y}$) \hfill $\lbrace$Orthogonalize column samples$\rbrace$}
		\STATE{$\bm{Z}=\bm{A}^{\top}\bm{Q_{Y}}$ \hfill $\lbrace$Sample row space$\rbrace$}
		\STATE{$\bm{Q_{Z}}=$ \verb~orthog~($\bm{Z}$) \hfill $\lbrace$Orthogonalize row samples$\rbrace$}
		\STATE{$\bm{R_{Z}}=\bm{Q_{Z}}^{\top}\bm{Z}$ \hfill $\lbrace$Compress row samples$\rbrace$}
		\STATE{$\bm{R_{Z}}=\bm{\hat{V}}\bm{\Sigma}\bm{\hat{U}}^{\top}$ \hfill $\lbrace$SVD of $k\times k$ compressed row sample  matrix$\rbrace$}
		\STATE{$\bm{V}=\bm{Q_{Z}}\bm{\hat{V}}$ \hfill $\lbrace$Project row space information$\rbrace$}
		\STATE{$\bm{U}=\bm{Q_{Y}}\bm{\hat{U}}$ \hfill $\lbrace$Project column space information$\rbrace$}
		\STATE{$\bm{\tilde{A}}=\bm{U}\bm{\Sigma}\bm{V}^{\top}$ \hfill $\lbrace$Form low-rank approximation$\rbrace$\,}
	\end{algorithmic}
	\label{alg:doublepass}
\end{algorithm}

The randomized hierarchical off-diagonal low-rank algorithm proceeds by compressing off-diagonal blocks by the double-pass algorithm. The larger off-diagonal blocks are compressed prior to the compression of smaller off-diagonal blocks, via a peeling procedure~\cite{linlulexing2011}. Here, both
$\bm{A}$ and $\bm{\tilde{A}}$ are assumed to be symmetric as we seek compression of symmetric operators and computation of symmetric approximants.
%
%
\begin{algorithm}[H]
	\caption{Symmetric matrix-free randomized HODLR. \\
		\textbf{Input:} symmetric $\bm{A}\in\mathbb{R}^{N\times N}$, hierarchical depth $L\in\mathbb{N}$, $r_{1},\dots,r_{L}$ desired ranks of the off-diagonal blocks at each hierarchical depth and oversampling parameter $d$.\\
		\textbf{Output:} symmetric HODLR approximation $\bm{\tilde{A}}$ of $\bm{A}$}
	\label{alg:HODLRcompression}
	\begin{algorithmic}[1]
		\FOR{$\ell=1,2,\dots,L$}
		\STATE{$k_{\ell}=r_{\ell}+d$}
		\STATE{$\bm{\Omega}=$ \verb~zeros~$(N,k_{\ell})$}
		\FOR{$j=1,\dots,2^{\ell-1}$}
		\STATE{$\bm{\Omega}(\mathcal{I}_{2j}^{(\ell)},:)=$ \verb~randn~$(|\mathcal{I}_{2j}^{(\ell)}|, k_{\ell})$ \hfill $\lbrace$Initiate structured random matrix$\rbrace$}
		\ENDFOR
		\STATE{$\bm{Y}=\left(\bm{A}-\sum_{j=1}^{\ell-1}\bm{A}^{(j)}\right)\bm{\Omega}$ \hfill $\lbrace$Sample off-diagonal block column spaces$\rbrace$}
		\FOR{$j=1,\dots,2^{\ell-1}$}
		\STATE{$\bm{Y}^{(j)}=$ zeros$(N,k_{\ell})$}
		\STATE{$\bm{Y}^{(j)}(\mathcal{I}_{2j-1}^{(\ell)},:)=\bm{Y}(\mathcal{I}_{2j-1}^{(\ell)},:)$}
		\STATE{$\bm{Q_{Y}}^{(j)}=$ \verb~orthog~$(\bm{Y}^{(j)})$ \hfill $\lbrace$Orthogonalize column samples of the level $\ell$ off-diagonal blocks$\rbrace$} 
		\ENDFOR
		\STATE{$\bm{Q_{Y}}=\sum_{j=1}^{2^{\ell-1}}\bm{Q_{Y}}^{(j)}$ \hfill $\lbrace$Row space sampling matrix$\rbrace$}
		\STATE{$\bm{Z}=\left(\bm{A}-\sum_{j=1}^{\ell-1}\bm{A}^{(j)}\right)\bm{Q_{Y}}$ \hfill $\lbrace$Sample off-diagonal block row spaces$\rbrace$}
		\FOR{$j=1,\dots,2^{\ell-1}$}
		\STATE{$\bm{Z}^{(j)}=\bm{Z}(\mathcal{I}_{2j}^{(\ell)},:)$}
		\STATE{$\bm{Q_{Z}}^{(j)}=$ \verb~orthog~$(\bm{Z}^{(j)})$ \hfill $\lbrace$Orthogonalize row samples of the level $\ell$ off-diagonal blocks$\rbrace$}
		\STATE{$\bm{R_{Z}}^{(j)}=\left(\bm{Q_{Z}}^{(j)}\right)^{\top}\bm{Z}^{(j)}$ \hfill $\lbrace$Compress level $\ell$ off-diagonal block row samples$\rbrace$}
		\STATE{$\bm{R_{Z}}^{(j)}=\bm{\hat{V}}_{2j-1}^{(\ell)}
			\bm{\Sigma}_{2j-1}^{(\ell)}
			\bm{\hat{U}}_{2j-1}^{(\ell)}$ \hfill $\lbrace$SVD of $k_{\ell}\times k_{\ell}$ compressed row sample matrix$\rbrace$ }
		\STATE{$\bm{V}_{2j-1}^{(\ell)}=\bm{Q_{Z}}^{(j)}\bm{\hat{V}}_{2j-1}^{(\ell)}$ \hfill $\lbrace$Project row space information$\rbrace$}
		\STATE{$\bm{U}_{2j-1}^{(\ell)}=\bm{Q_{Y}}^{(j)}\bm{\hat{U}}_{2j-1}^{(\ell)}$ \hfill $\lbrace$Project column space information$\rbrace$}
		\STATE{$\bm{V}_{2j}^{(\ell)}=\bm{U}_{2j-1}^{(\ell)}$}
		\STATE{$\bm{U}_{2j}^{(\ell)}=\bm{V}_{2j-1}^{(\ell)}$}
		\STATE{$\bm{\Sigma}_{2j}^{(\ell)}=\bm{\Sigma}_{2j-1}^{(\ell)}$}
		\ENDFOR
		\STATE{$\bm{A}^{(\ell)}=\sum_{j=1}^{2^{\ell}}\bm{U}_{j}^{(\ell)}\bm{\Sigma}_{j}^{(\ell)}
			\left(\bm{V}_{j}^{(\ell)}\right)^{\top}$} 		 
		\ENDFOR
		\STATE{obtain block diagonal $\bm{D}$ of $\bm{A}$ by sampling $\bm{A}-\sum_{j=1}^{L}\bm{A}^{(j)}$}
		\STATE{$\bm{\tilde{A}}=\bm{D}+\sum_{\ell=1}^{L}\bm{A}^{(\ell)}$}
	\end{algorithmic}
\end{algorithm}

\subsection{Global HODLR approximation error from the accumulation of block low-rank off-diagonal approximation errors}
\label{subsec:erroranalysis}

Let $\bm{A}$ be a $N\times N$ matrix and consider the following partitioning
\begin{eqnarray*}
\bm{A}^{\left(1\right)}
&=&
\left(\matrix{
\bm{0} & \bm{A}_{1,2}^{\left(1\right)} \cr
\bm{A}_{2,1}^{\left(1\right)} & \bm{0} \cr}\right), \\
\bm{A}^{\left(2\right)}
&=&
\left(\matrix{
\bm{0} & \bm{A}_{1,2}^{\left(2\right)} &
\bm{0} & \bm{0} \cr
\bm{A}_{2,1}^{\left(2\right)} & \bm{0} &
\bm{0} & \bm{0} \cr
\bm{0} & \bm{0} & 
\bm{0} & \bm{A}_{3,4}^{\left(2\right)} \cr
\bm{0} & \bm{0} & 
\bm{A}_{4,3}^{\left(2\right)} & \bm{0} \cr
}\right), \\
\bm{D}&=&
\left(\matrix{
\bm{A}_{1,1}^{\left(2\right)} & \bm{0} &
\bm{0} & \bm{0} \cr 
\bm{0} & \bm{A}_{2,2}^{\left(2\right)} &
\bm{0} & \bm{0} \cr
\bm{0} & \bm{0} &
\bm{A}_{3,3}^{\left(2\right)} & \bm{0} \cr
\bm{0} & \bm{0} &
\bm{0} & \bm{A}_{4,4}^{\left(2\right)} \cr
}\right),
\end{eqnarray*}
where $\bm{A}_{i,j}^{\left(\ell\right)}$ is the $(i,j)$ block of a $2^{\ell}\times 2^{\ell}$ 
block partitioning of $\bm{A}$, where $1\leq \ell\leq L$. $\bm{A}^{(\ell)}$ contains all blocks $\bm{A}_{i,j}^{(\ell)}$ such that $|i-j|=1$ and $\bm{D}$ contains the 
diagonal blocks $\bm{A}^{(L)}_{i,i}$. Above, we show the decomposition $\bm{A}=\sum_{\ell=1}^{L}\bm{A}^{(\ell)}+\bm{D}$ for $L=2$ hierarchical depth but in the following analysis 
$L$ is a arbitrary. Let $\bm{x}\in\mathbb{R}^{N}$, then
\begin{eqnarray*}
\bm{A}\bm{x}&=&
\sum_{j=1}^{L}\bm{A}^{\left(j\right)}\bm{x}+\bm{D}\bm{x}, \\
\bm{A}^{\left(1\right)}\bm{x}
&=&
\left(\matrix{
\bm{A}_{1,2}^{\left(1\right)}\bm{x}_{2}^{\left(1\right)} \cr
\bm{A}_{2,1}^{\left(1\right)}\bm{x}_{1}^{\left(1\right)} \cr 
}\right), 
\,\,\,\bm{x}=
\left(\matrix{
\bm{x}_{1}^{\left(1\right)} \cr 
\bm{x}_{2}^{\left(2\right)} \cr 
}\right), \\
\bm{A}^{\left(j\right)}\bm{x}&=&
\left(\matrix{
\bm{A}_{1,2}^{\left(j\right)}\bm{x}_{2}^{\left(j\right)} \cr 
\bm{A}_{2,1}^{\left(j\right)}\bm{x}_{1}^{\left(j\right)} \cr 
\vdots \cr 
\bm{A}_{2^{j}-1,2^{j}}^{\left(j\right)}\bm{x}_{2^{j}}^{\left(j\right)} \cr
\bm{A}_{2^{j},2^{j}-1}^{\left(j\right)}\bm{x}_{2^{j}-1} \cr
}\right), 
\,\,\,\bm{x}
=
\left(\matrix{
\bm{x}_{1}^{\left(j\right)} \cr 
\bm{x}_{2}^{\left(j\right)} \cr
\vdots \cr
\bm{x}_{2^{j}-1}^{\left(j\right)} \cr
\bm{x}_{2^{j}}^{\left(j\right)} \cr
}\right),
\end{eqnarray*}
from which we obtain the following expression
\begin{equation*}
\|\bm{A}^{\left(j\right)}\bm{x}\|_{2}^{2}=
\sum_{k=1}^{2^{j-1}}\left(
\|\bm{A}_{2\,k-1,2\,k}^{\left(j\right)}
\bm{x}_{2\,k}^{\left(j\right)}\|_{2}^{2}
+
\|\bm{A}_{2\,k,2\,k-1}^{\left(j\right)}
\bm{x}_{2\,k-1}^{\left(j\right)}\|_{2}^{2}
\right). 
\end{equation*}
Now assume that $\bm{\tilde{A}}$ is an HODLR approximation of $\bm{A}$, whose diagonal $\bm{D}$ is equal to the diagonal of $\bm{A}$ so that
\begin{eqnarray*}
\left(\bm{A}-\bm{\tilde{A}}\right)
&=&\sum_{j=1}^{L}
\Delta\bm{A}^{\left(j\right)},\\
\Delta\bm{A}^{\left(j\right)}&:=&\left(\bm{A}^{\left(j\right)}-\bm{\tilde{A}}^{\left(j\right)}\right).
\end{eqnarray*}
Here we assume each off-diagonal block has been approximated to some absolute tolerance $\varepsilon>0$,
so that $\|\Delta \bm{A}_{2\,k-1,2\,k}^{\left(j\right)}\|_{2},
\|\Delta \bm{A}_{2\,k,2\,k-1}^{\left(j\right)}\|\leq\varepsilon$ for each $j= 1,2,\dots, L$ and $k = 1,2,\dots,2^{j-1}$. For $\bm{x}\in\mathbb{R}^{N}$ we have
\begin{equation*}
\|\left(\bm{A}-\bm{\tilde{A}}\right)\bm{x}\|_{2}
\leq 
\sum_{j=1}^{L}\|\Delta\bm{A}^{\left(j\right)}\,\bm{x}\|_{2},
\end{equation*}
\begin{eqnarray*}
\|\Delta\bm{A}^{\left(j\right)}\,\bm{x}\|_{2}&=&
\sqrt{ 
	\sum_{k=1}^{2^{j-1}}\left(
	\|\Delta \bm{A}_{2\,k-1,2\,k}^{\left(j\right)}
	\,\bm{x}_{2\,k}^{\left(j\right)}\|_{2}^{2}
	+
	\|
	\Delta 
	\bm{A}_{2\,k,2\,k-1}^{\left(j\right)}\,
	\bm{x}_{2\,k-1}^{\left(j\right)}\|_{2}^{2}
	\right) 
} \\
&\leq&
\sqrt{
	\sum_{k=1}^{2^{j-1}}
	\left(\varepsilon^{2}
	\|\bm{x}_{2\,k}^{\left(j\right)}\|_{2}^{2}
	+\varepsilon^{2}
	\|\bm{x}_{2\,k-1}^{\left(j\right)}\|_{2}^{2}\right) 
}, \\
\|\Delta \bm{A}^{\left(j\right)}\,\bm{x}\|_{2}
&\leq&
\varepsilon 
\sqrt{
	\sum_{k=1}^{2^{j-1}}
	\left(
	\|\bm{x}_{2\,k}^{\left(j\right)}\|_{2}^{2}
	+\|\bm{x}_{2\,k-1}^{\left(j\right)}\|_{2}^{2}\right)
}
=\varepsilon\|\bm{x}\|_{2}, \\
\|\left(\bm{A}
-\bm{\tilde{A}}\right) 
\bm{x}\|_{2}
&\leq& \varepsilon \,L\, \|\bm{x}\|_{2}, \\
\|\bm{A}-\bm{\tilde{A}}\|_{2}&:=&
\sup_{\bm{x}\neq\bm{0}}
\left( 
\frac{\|\left(\bm{A}-\bm{\tilde{A}}\right)\bm{x}\|_{2}}
{\|\bm{x}\|_{2}}\right)\leq \varepsilon\,L .
\end{eqnarray*}
\subsection{Error analysis for posterior-covariance}
\label{subsec:posteriorerroranalysis}
Consider a symmetric matrix 
$\bm{A}\in\mathbb{R}^{N\times N}$, whose eigenvalues are bounded below by a number greater than $-1$ and a symmetric 
approximant $\bm{\tilde{A}}$, with discrepancy $\Delta\bm{A}=\bm{A}-\bm{\tilde{A}}$.
We signify a generic eigenvalue of $\bm{S}$ by $\lambda\left(\bm{S}\right)$ so that $s_{1}\leq \lambda\left(\bm{S}\right)\leq s_{2}$ indicates that all eigenvalues of $\bm{S}$ are bounded below by $s_{1}$ and above by $s_{2}$.  
Now we provide a bound for the error of $\left(\bm{I}+\bm{A}\right)^{-1}-\left(\bm{I}+\bm{\tilde{A}}\right)^{-1}$,
given that $\|\Delta\bm{A}\|_{2}= \varepsilon$, so that one 
may assess the accuracy of an HODLR Gaussianized posterior covariance.
When, as in Section~\ref{subsec:HODLRGaussianizedPosterior}, $\bm{A}$ is the prior-preconditioned Hessian misfit, $\|\left(\bm{I}+\bm{A}\right)^{-1}-\left(\bm{I}+\bm{\tilde{A}}\right)^{-1}\|_{2}$
quantifies the discrepancy between an HODLR approximate Gaussianized posterior covariance and the true Gaussianized posterior covariance. 
\begin{eqnarray*}
\left(\bm{I}+\bm{A}\right)^{-1}-\left(\bm{I}+\bm{\tilde{A}}\right)^{-1}=\left(\bm{I}+\bm{A}\right)^{-1}-\left(\bm{I}+\bm{A}-\Delta\bm{A}\right)^{-1} =\\
\left(\bm{I}+\bm{A}\right)^{-1}-\left(\left(\bm{I}+\bm{A}\right)\left(\bm{I}-
\left(\bm{I}+\bm{A}\right)^{-1}
\Delta\bm{A}\right)\right)^{-1}=\\
\left(\bm{I}+\bm{A}\right)^{-1}-\left(\bm{I}-
\left(\bm{I}+\bm{A}\right)^{-1}
\Delta\bm{A}\right)^{-1}\left(\bm{I}+\bm{A}\right)^{-1}=\\
\left(\bm{I} - \left(\bm{I}-
\left(\bm{I}+\bm{A}\right)^{-1}
\Delta\bm{A}\right)^{-1}\right)\left(\bm{I}+\bm{A}\right)^{-1}.
\end{eqnarray*}
Given that $\|\Delta\bm{A}\|_{2}=\varepsilon$, we have
\begin{eqnarray*}
-\varepsilon \leq \lambda\left(\Delta\bm{A}\right)\leq \varepsilon, \\
-\varepsilon^{*}
\leq \lambda\left(\left(\bm{I}+\bm{A}\right)^{-1}\Delta\bm{A}\right)\leq \varepsilon^{*}, \\
\varepsilon^{*}:=\varepsilon(1+\lambda_{\text{min}}(\bm{A}))^{-1},\\
1+\varepsilon^{*} \geq \lambda\left(\bm{I}-\left(\bm{I}+\bm{A}\right)^{-1}\Delta\bm{A}\right)\geq 1- \varepsilon^{*},
\end{eqnarray*}
we next assume $\varepsilon^{*}<1$, so that the eigenvalues of $\bm{I}-\left(\bm{I}+\bm{A}\right)^{-1}\Delta\bm{A}$ are necessarily positive and
\begin{equation*}
\left(1+\varepsilon^{*}\right)^{-1}
\leq \lambda\left(\left(\bm{I}-\left(\bm{I}+\bm{A}\right)^{-1}\Delta\bm{A}\right)^{-1}\right)\leq \left(1-\varepsilon^{*}\right)^{-1}.
\end{equation*}
With this it follows that
\begin{eqnarray*}
\|\left(\bm{I}+\bm{A}\right)^{-1}-
\left(\bm{I}+\bm{\tilde{A}}\right)^{-1}\|_{2}/\|\left(\bm{I}+\bm{A}\right)^{-1}\|_{2}\leq \left(1-\left(1+\varepsilon^{*}\right)^{-1}\right) \\
\|\left(\bm{I}+\bm{A}\right)^{-1}-
\left(\bm{I}+\bm{\tilde{A}}\right)^{-1}\|_{2}/\|\left(\bm{I}+\bm{A}\right)^{-1}\|_{2}\leq \frac{\varepsilon^{*}}{1+\varepsilon^{*}},
\end{eqnarray*}
where, as before $\varepsilon^{*}=\|\Delta\bm{A}\|_{2}/\left(1+\lambda_{\text{min}}\left(\bm{A}\right)\right)$.

\section*{Acknowledgments}
The authors thank Trevor Hillebrand from Los Alamos National Laboratory for help with setting up the Humboldt and Greenland ice-sheet grids and datasets. Support for this work was provided by the National Science Foundation under Grant No.\ DMS-1840265 and CAREER-1654311 and through the SciDAC project ProSPect, funded by the U.S.\ Department of Energy (DOE) Office of Science, Advanced Scientific Computing Research and Biological and Environmental Research programs. This research used resources of the National Energy Research Scientific Computing Center (NERSC), a U.S. Department of Energy Office of Science User Facility operated under Contract No.\ DE-AC02-05CH11231, under NERSC award ERCAP0020130.

\section*{Disclaimer}
This paper describes objective technical results and analysis. Any subjective views or opinions that might be expressed in the paper do not necessarily represent the views of the  U.S. Department of Energy or the United States Government. Sandia National Laboratories is a multimission laboratory managed and operated by National Technology and Engineering Solutions of Sandia, LLC, a wholly owned subsidiary of Honeywell International, Inc., for the U.S. Department of Energy’s National Nuclear Security Administration under contract DE-NA-0003525.

\section*{References}

\bibliographystyle{iopart-num}
\bibliography{references}
\end{document}